\documentclass[12pt,draft]{article}

\usepackage{amsthm,amssymb,amsmath,amscd,amstext,amsopn,mathrsfs}
\usepackage{latexsym,mathptmx,xypic,verbatim}

\newtheorem{theorem}{Theorem}[section]

\newtheorem{lemma}{Lemma}[section]
\newtheorem*{remark}{{\it Remark}}

\DeclareMathOperator{\tr}{tr}
\DeclareMathOperator{\Ad}{Ad}
\DeclareMathOperator{\ad}{ad}
\DeclareMathOperator{\re}{Re}
\DeclareMathOperator{\im}{Im}
\DeclareMathOperator{\id}{Id}

\newcommand{\nc}{\newcommand}

\nc{\OO}{{\mathbb O}}
\nc{\C}{{\mathbb C}}
\nc{\R}{{\mathbb R}}
\nc{\Z}{{\mathbb Z}}
\nc{\N}{{\mathbb N}}

\nc{\ee}{{\bf e}}
\nc{\dd}{{\rm d}}

\begin{document}

\title{The complex structure on the six dimensional sphere}

\author{G\'abor Etesi\\
\small{{\it Department of Algebra and Geometry, Institute of Mathematics,}}\\
\small{{\it Budapest University of Technology and Economics,}}\\
\small{{\it M\H uegyetem rkp. 3., H-1111 Budapest, Hungary}}
\footnote{E-mail: {\tt etesi@math.bme.hu, etesigabor@gmail.com}}}
\maketitle

\pagestyle{myheadings}
\markright{G. Etesi: The complex structure on the six-sphere}

\thispagestyle{empty}

\begin{abstract} 
Proof of existence of at least one complex structure on 
the six-sphere, followed by an explicit computation of its underlying 
integrable almost complex tensor by the aid of inner automorphisms of the 
octonions, is exhibited. Both are elementary and self-contained however the 
size and complexity of the emerging almost complex tensor field on the 
six-sphere is perplexing. 
\end{abstract}

\centerline{AMS Classification: Primary: 32J17, 32J25; Secondary: 22C05, 53C56}
\centerline{Keywords: {\it Six-sphere; Complex structure; Octonions}}


\section{Introduction}
\label{one}


Unlike their real counterparts, complex manifolds are quite rigid objects 
therefore constructing new examples is often not an easy task. Nevertheless 
there are two main islands, both within the {\it K\"ahler realm}, of the 
archipelago of complex manifolds comprising tractable cases. One is the 
category of {\it projective manifolds}, standing in the focal point of 
algebraic geometry. The other one is the category of {\it Stein manifolds} 
which, on the contrary, can be conveniently studied with the techniques of 
complex analysis in several variables. Although the in-between 
{\it terra incognita} is not so easily accessable, it is still quite 
populated: its habitants are the numerous irregular neither projective nor 
Stein---or even not K\"ahler---manifolds, either compact or non-compact. 
Examples in two complex dimensions are for example the various non-algebraic 
tori, Hopf surfaces, Inue surfaces, non-algebraic $K3$ surfaces, etc. 
(cf., e.g. \cite[Chapter VI]{bar-hul-pet-ven} for a survey of irregular 
complex surfaces) while in higher dimensions the picture is not so clear yet. 
Nevertheless in all dimensions the deviation, at least from the algebraic 
scenario, is captured in some extent by the concept of the algebraic dimension 
of a complex manifold i.e. the transcendental degree over $\C$ of its field of 
global meromorphic functions. Compact complex manifolds whose algebraic and 
geometric dimensions match (sometimes called {\it Moishezon manifolds}) 
are still not far from being algebraic: it is known that after performing 
finitely many blowing-ups they become projective algebraic. However in 
general the lower the algebraic dimension of a manifold is the larger 
its detachment is from the familiar world of complex algebraic manifolds. 

An interesting family of higher dimensional non-algebraic examples 
relevant to us here is based on even dimensional compact Lie groups. 
{\it Samelson} discovered \cite{sam} that they carry complex structures 
which are surely non-K\"ahler if the underlying group is simply connected and 
simple. It is therefore guessed that in general they admit lower algebraic 
dimension than possible. It also has been known 
for some time that if the six dimensional sphere carries a complex structure 
then this compact complex $3$-manifold must be of zero algebraic dimension 
\cite{cam-dem-pet}. 

In this paper we shall construct explicitly at least one complex structure 
on the six dimensional sphere. The proof is based on identifying $S^6$ with 
an exceptional conjugate orbit in the exceptional compact Lie group 
${\rm G}_2$, taking its explicit deformation within 
${\rm G}_2^\C$ and then restricting a Samelson complex structure to this 
deformation. The proof 
to be presented here is elementary and self-contained hence is 
independent of our former Yang--Mills--Higgs theoretic approach 
\cite{ete}; nevertheless those considerations definitely have been used 
here as a source of ideas. We just note that meanwhile the treatment in 
\cite{ete} is based on the well-known ${\rm SU}(3)$-fibration: the {\it 
projection} (i.e. a surjective mapping) $\pi: {\rm G}_2\rightarrow S^6$, 
our present proof rests on a certain less-known but very remarkable {\it 
injection} $f:S^6\rightarrow{\rm G}_2$. 
Therefore the two approaches are dual to each other in this sense.

The paper is organized as follows. Sect. \ref{two} consists of a 
self-contained and elementary proof that $S^6$ can carry complex 
structures (cf. Theorem \ref{letezes}). It has been known for a long 
time (cf. e.g. \cite{leb, huc-keb-pet}) that if $S^6$ carries a complex 
structure then there exist ``exotic $\C P^3$'s'' i.e. complex manifolds 
diffeomoprhic to $\C P^3$ as real $6$-manifolds however not 
complex-analytically isomorphic to it. Here we prove another striking 
consequence: the complex sructure on $S^6$ implies the existence of 
``large exotic $\C^3$'s'' in a similar sense (cf. Lemma 
\ref{egzotikus}). By a {\it large exotic $\C^3$} we mean a complex 
manifold which is diffeomorphic to $\R^6$ however is not complex 
analytically equivalent to the standard $\C^3$ moreover does not admit a 
complex-analytic embedding into the standard $\C^3$ (the failure of the 
higher dimensional analogue of the Riemannian mapping theorem implies 
the existence of an abundance of {\it small exotic $\C^3\:$'s} i.e. open 
complex analytic subsets of the standard $\C^3$ not complex-analytically 
equivalent to it).

In Sect. \ref{three} we outline the explicit construction of the 
integrable almost complex tensor field on $S^6$ underlying its complex 
structure. The construction is based on identifying the original 
conjugate orbit of ${\rm G}_2$, homeomorphic to $S^6$, with the subset of 
{\it inner automorphisms} within the full automorphism group of the octonions 
which is ${\rm G}_2$ (cf. Lemma \ref{egyenloseg} as well as \cite{cha-rig}) 
and then taking its explicit perturbation within ${\rm G}_2^\C$ 
(cf. Lemma \ref{komplexegyenloseg}). The integrable Samelson 
almost complex structures then restrict to this perturbation rendering $S^6$ 
a complex manifold. The components of the underlying almost 
complex tensor field considered as local $6\times 6$ matrix functions 
in principle drop out explicitly from this construction 
however the result is so unexpectedly complicated that we cannot display it 
fully here. Nevertheless the steps towards its construction are clearly 
explained and the curious reader can reproduce the calculations (using a 
computer is strongly advised) by himself or even go further and bring these 
matrices into a more digestable form. But already the present form 
allows to conclude that all complex structures constructed here are equivalent.

Finally Sect. \ref{four} is an appendix and has been added in order to 
gain a more comprehensive picture. Following \cite{cha-rig} we re-prove 
that the conjugate orbit of ${\rm G}_2$ playing the central role here, when 
regarded as a continuous map $f:S^6\rightarrow{\rm G}_2$, represents the 
generator of $\pi_6({\rm G}_2)\cong\Z_3$ (cf. Theorem \ref{homotopia}).
\vspace{0.1in}

\noindent{\bf Acknowledgement.} The author is grateful to F. Burstall, 
B. Csik\'os, Sz. Szab\'o and R. Sz\H oke for the criticism of earlier versions 
and the stimulating discussions many years ago. The math software Maple 
has been extensively used to carry out the massive but {\it 
strictly symbolic} calculations in Sect. \ref{three}.


\section{Proof of existence}
\label{two}


In this section we present a proof that the six dimensional 
sphere carries complex structures. The proof is based on identifying $S^6$ with 
a geometrically deformed conjugate orbit inside the complexified exceptional 
compact Lie group ${\rm G}_2^\C$ and then restricting a Samelson complex 
structure to this deformed orbit. The proof to be 
presented here is elementary and self-contained hence is independent of our 
former Yang--Mills--Higgs theoretic approach \cite{ete}.  

Recall \cite{pit, sam} that if $G$ is an even dimensional compact real 
Lie group with real Lie algebra ${\mathfrak g}$ the complex linear 
subspace ${\mathfrak s}\subset{\mathfrak g}^\C:={\mathfrak 
g}\otimes_\R\C$ is a {\it Samelson subalgebra} if it is a complex Lie 
subalgebra of ${\mathfrak g}^\C$ satisfying $\dim_\C{\mathfrak 
s}=\frac{1}{2}\dim_\C{\mathfrak g}^\C$ and ${\mathfrak s}\cap{\mathfrak 
g}=0$ as real subspaces within ${\mathfrak g}^\C$. One can demonstrate 
that at least one Samelson subalgebra in ${\mathfrak g}^\C$ always 
exists if ${\mathfrak g}$ is the Lie algebra of a {\it compact} even 
dimensional Lie group. A choice for a Samelson subalgebra gives rise to 
a vector space decomposition ${\mathfrak g}^\C={\mathfrak 
s}\oplus\overline{{\mathfrak s}}$ hence the existence of a real linear 
isomorphism $\re: {\mathfrak s} \rightarrow{\mathfrak g}$ given by 
$W\mapsto\re W$ for all $W\in{\mathfrak s}$ and a real linear map 
$J_{\mathfrak s}:{\mathfrak g}\rightarrow{\mathfrak g}$ with 
$J_{\mathfrak s}(\re W):=-\im W$ satisfying $J_{\mathfrak 
s}^2=-\id_{\mathfrak g}$. Consequently a choice of a Samelson subalgebra 
gives rise to a complex vector space $({\mathfrak g},J_{\mathfrak s})$ 
such that ${\mathfrak g}^{1,0}:= {\mathfrak s}$ is the 
$+\sqrt{-1}\:$-eigenspace while ${\mathfrak g}^{0,1}:= 
\overline{{\mathfrak s}}$ is the $-\sqrt{-1}\:$-eigenspace of the 
complexified map $J_{\mathfrak s}^\C:{\mathfrak g}^\C\rightarrow 
{\mathfrak g}^\C$. The operator $J_{\mathfrak s} =J_{{\mathfrak s},e}$ 
on ${\mathfrak g}=T_eG$ can be extended to a left-invariant almost 
complex structure $J_{\mathfrak s}$ over the whole $G$ if $J_{{\mathfrak 
s},g}$ on $T_gG$ is defined by left-translating everything i.e. by the 
aid of the complexified group $G^\C$ taking the splitting 
$T_gG^\C=L_{g*}{\mathfrak s}\oplus L_{g*}\overline{\mathfrak s}$ and 
observing that $T_gG^\C=T_gG\otimes_\R\C$ such that $L_{g*}{\mathfrak 
s}\cap T_gG=0$ inside $T_gG^\C$ for every $g\in G\subset G^\C$. Note 
that simply $J_{{\mathfrak s},g}=L_{g\:*}\:J_{{\mathfrak s},e}\:L^{-1}_{g\:*}$. 
In this way we come up 
with an almost complex manifold $(G,J_{\mathfrak s})$. The Lie algebra 
property of ${\mathfrak g}^{1,0}$ has not been exploited so far: Taking 
into account ${\mathfrak g}^\C ={\mathfrak g}^{1,0}\oplus{\mathfrak 
g}^{0,1}$ as well, it additionally tells us that \[\left[ {\mathfrak 
g}^{1,0}\:,\:{\mathfrak g}^{1,0}\right]^{0,1}=0\] also holds. The 
identification $X\mapsto X^{1,0}:=\frac{1}{2}(X-\sqrt{-1}\:J_{\mathfrak 
s}X)$ of real vector fields with $(1,0)$-type complex ones over $G$ maps 
left-invariant real fields into left-invariant $(1,0)$-type ones and 
these latter fields can be viewed as $(1,0)$-type complex Lie algebra 
elements. Observing that $(X,Y)\mapsto [X^{1,0}, Y^{1,0}]^{0,1}$ is 
$C^\infty (G;\R )$-bilinear we recognize that the vanishing commutator 
above actually says that the Nijenhuis tensor of $J_{\mathfrak s}$ is 
zero. Consequently $(G, J_{\mathfrak s})$ is integrable to a homogeneous 
complex manifold $Y_{\mathfrak s}$ in light of the Newlander--Nirenberg 
theorem. (For a different proof cf. \cite[Proposition 2.3]{pit}.) 

After these general considerations take the $14$ dimensional real 
compact exceptional Lie group ${\rm G}_2$. As we outlined above it can be given 
the structure of a compact complex $7$-manifold; more precisely referring to 
\cite[Example on p. 123]{pit} we know\footnote{But we will also 
recover this result explicitly in Sect. \ref{three}.} that if ${\mathfrak h}
\subset{\mathfrak g}_2$ is the Cartan subalgebra with its 
complexification ${\mathfrak h}^\C\subset{\mathfrak g}_2^\C$ then for 
every $u\in P({\mathfrak h}^\C)\setminus P({\mathfrak h})$ i.e. the 
projectivization of ${\mathfrak h}^\C$ with that of its real part removed, 
there exist Samelson subalgebras ${\mathfrak s}_u\subset{\mathfrak g}_2^\C$ 
with corresponding mutually non-isomorphic compact homogeneous complex 
$7$-manifolds $Y_u$ whose underlying real spaces are all diffeomorphic to 
${\rm G}_2$. This Lie group has a maximal subgroup isomorphic to 
${\rm SU}(3)$ and let $\Lambda\in {\rm G}_2$ be the generator of its center 
i.e. $\Lambda\in Z({\rm SU}(3))\subset{\rm SU}(3)\subset{\rm G}_2$ is 
the generator, hence $\Lambda\not=e$ but $\Lambda^3=e$. Strongly motivated by 
\cite{cha-rig, dur-put-rig} we begin with considering the exceptional 
{\it conjugate orbit} 
\begin{equation}
O(\Lambda ):=\left\{ g\Lambda g^{-1}\:\vert\:g\in {\rm G}_2\right\}
\label{palya}
\end{equation}
passing through this generator. It is a real submanifold of ${\rm G}_2$ 
diffeomorphic to ${\rm G}_2/{\rm SU}(3)\cong S^6$. Let ${\rm G}_2^\C$ be the 
complexification of ${\rm G}_2$. Our real conjugate orbit
$O(\Lambda )\subset {\rm G}_2$ also complexifies in an obvious way to
the complexified conjugate orbit
\[O(\Lambda)^\C=\{g\Lambda g^{-1}\:\vert\:g\in{\rm G}_2^\C\}\:\:.\]
By compactness of ${\rm G}_2$ there exists a diffeomorphism
${\rm G}_2^\C\cong T{\rm G}_2$ as a real manifold hence we can use the zero
section of the tangent bundle to write ${\rm G}_2\subset{\rm G}_2^\C$. This 
by restriction gives $O(\Lambda)\subset O(\Lambda)^\C$ and the 
real isomorphism above implies $O(\Lambda )^\C\cong TO(\Lambda )$.
Note that $O(\Lambda)^\C$ is complex-analytically isomorphic to
${\rm G}_2^\C/{\rm SU}(3)^\C\cong (S^6)^\C$ where $(S^6)^\C\subset\C^7$ is
the ``complex 6-sphere'' whose points satisfy $z_1^2+\dots +z_7^2=1$ with
$z_1, \dots,z_7\in\C$. Therefore our embeddings are in accordance 
with the classical fact that $(S^6)^\C\cong TS^6$ as a real manifold and 
with the existence of an embedding $S^6\subset (S^6)^\C$ by the zero section. 
To summarize, there exists a commutative diagram 
\begin{equation}
\begin{matrix} 
{\rm G}_2&\subset &{\rm G}_2^\C&\cong&T{\rm G}_2&\cong&
                                       {\rm G}_2\times{\mathfrak g}_2\\
                  \cup    &        &    \cup    &     &  \cup    & &  \\
      O(\Lambda ) & \subset & O(\Lambda )^\C &\cong & TO(\Lambda ) & &\\    
          \wr\!\!\parallel & & \wr\!\!\parallel & &\wr\!\!\parallel & &\\
                         S^6    & \subset&  (S^6)^\C   &\cong& TS^6 & &
\end{matrix}
\label{lekepezesek}
\end{equation}
of compatible embeddings and isomorphisms. 

Next we perform a complex deformation of $O(\Lambda )$ within ${\rm G}^\C_2$. 
Take a fixed Samelson subalgebra ${\mathfrak s}_u\subset{\mathfrak g}_2^\C$ 
with corresponding complex analytic Lie subgroup $S_u\subset {\rm G}_2^\C$ and 
likewise take their complex conjugate counterparts 
$\overline{\mathfrak s}_u\subset{\mathfrak g}_2^\C$ and 
$\overline{S}_u\subset{\rm G}_2^\C$. Consider a new complex analytic 
submanifold
\begin{equation}
O'(\Lambda )^\C:=\left\{s_1\Lambda^2s_1^{-1}s_2\Lambda^2 
s_2^{-1}\:\big\vert\:s_1\in S_u\:,\: s_2\in \overline{S}_u\right\}
\label{komplexpalya'}
\end{equation}
inside ${\rm G}_2^\C$ and its ``real part'' $O'(\Lambda )\subset 
O'(\Lambda )^\C$ defined by
\begin{equation}
O'(\Lambda ):=\left\{ s\Lambda^2s^{-1}\:\overline{s}\Lambda^2
\overline{s}^{-1}\:\big\vert\:s\in S_u\right\}
\label{palya'}
\end{equation} 
and regarded as a {\it deformed conjugate orbit} through $\Lambda$ (note that 
$\Lambda\in O'(\Lambda )$ continues to hold) but 
within ${\rm G}_2^\C$ (note that $O'(\Lambda )\not\subset{\rm G}_2$). 
We assert that $O'(\Lambda )$ is homeomorphic to $S^6$. 

To prove this first we observe that 
$O(\Lambda )^\C$ and $O'(\Lambda )^\C$ are isotopic in ${\rm G}_2^\C$. 
The decomposition ${\mathfrak g}_2^\C =
{\mathfrak s}_u\oplus\overline{\mathfrak s}_u$ implies ${\rm G}_2^\C
\supset S_u\overline{S}_u$ is open and 
${\rm G}_2^\C\supset S_u\cap \overline{S}_u$ is discrete; thus 
$S_u\Lambda^2S_u^{-1}\cap\overline{S}_u\Lambda^2\overline{S}_u^{-1}\subset
{\rm G}_2^\C$ is discrete and since both $S_u\cap{\rm SU}(3)^\C$ and
$\overline{S}_u\cap {\rm SU}(3)^\C$ are $4$ dimensional and $\Lambda^2$
commutes with them $\dim_\C O'(\Lambda )^\C = 2\cdot (7-4)=6$. Take 
the $1$-parameter subgroups $\sigma_1:\R\rightarrow S_u$ and 
$\sigma_2:\R\rightarrow\overline{S}_u$ satisfying $\sigma_i(0)=e$ and 
$\sigma_i(1)=s_i$. Observe that $S_u$ is solvable \cite[Proposition 2.4]{pit} 
and simply connected hence $\exp:{\mathfrak s}_u\rightarrow S_u$ 
is a diffeomorphism 
\cite[Theorem 3.18.11]{var} and likewise for $\overline{S}_u$. Consequently 
$\sigma_i$ uniquely exist thus for every $n\in\N$ and $s_i$ one can 
unambigously put $s_i^{1/n}:=\sigma_i(\frac{1}{n})$ and hence define 
$S^{1/n}_u:=\{s_1^{1/n}\:\vert\:s_1\in S_u\}$ and likewise for 
$\overline{S}_u^{1/n}$. Furthermore ${\rm G}_2^\C$ is a complex 
simple group with trivial center hence $\exp:{\mathfrak g}_2^\C\rightarrow 
{\rm G}_2^\C$ is surjective \cite[Corollary 3.4]{mos} but not injective. 
Consequently, recalling Trotter's formula $\exp(X_1+X_2)=
\lim\limits_n\big(\exp\frac{X_1}{n}\exp\frac{X_2}{n}\big)^n$ from e.g. 
\cite[Corollary 2.12.5]{var}, any $g\in {\rm G}_2^\C$ can be written
as $g=\lim\limits_n\big(s_1^{1/n}s_2^{1/n}\big)^n$ that is, 
${\mathfrak g}_2^\C={\mathfrak s}_u\oplus\overline{\mathfrak s}_u$ globalizes 
to ${\rm G}_2^\C=\lim\limits_n\big(S^{1/n}_u
\:\:\overline{S}_u^{1/n}\big)^n$ but in a non-unique way. 
However uniqueness is achieved for the quotient ${\rm G}_2^\C/{\rm SU}(3)^\C$. 
Indeed, take a small open subset $U\subset{\mathfrak g}_2^\C$ about the 
origin and consider the local model 
$U/{\mathfrak s}{\mathfrak u}(3)^\C$ for ${\rm G}_2^\C/{\rm SU}(3)^\C$. 
Picking elements $a\in\exp U\subset{\rm G}_2^\C$ and 
$X,Y\in\exp^{-1}(a{\rm SU}(3)^\C)$, since ${\rm SU}(3)^\C\subset{\rm G}_2^\C$ 
has maximal rank hence ${\mathfrak s}{\mathfrak u}(3)^\C$ contains the real 
integer lattice of ${\mathfrak g}_2^\C$ which generates the 
kernel of the relative exponential map 
$\exp:({\mathfrak g}_2^\C,{\mathfrak s}{\mathfrak u}(3)^\C)
\rightarrow({\rm G}_2^\C,{\rm SU}(3)^\C)$ we find that
$X+{\mathfrak s}{\mathfrak u}(3)^\C=Y+{\mathfrak s}
{\mathfrak u}(3)^\C$ i.e. $Y\in X+{\mathfrak s}{\mathfrak u}(3)^\C$; since 
$\exp U$ generates ${\rm G}_2^\C$ we conclude that 
$X,Y\in\exp^{-1}(g{\rm SU}(3)^\C)$ implies 
$Y\in X+{\mathfrak s}{\mathfrak u}(3)^\C$ for $g\in{\rm G}_2^\C$. 
Thus if $X=X_1+X_2$ is the unique Samelson decomposition then surely 
$Y=(X_1+V_1)+(X_2+V_2)$ where $V_1\in{\mathfrak s}_u\cap
{\mathfrak s}{\mathfrak u}(3)^\C$ and 
$V_2\in\overline{\mathfrak s}_u\cap{\mathfrak s}{\mathfrak u}(3)^\C$ hence 
using Trotter's formula we obtain that to {\it any} coset 
$g{\rm SU}(3)^\C$ one can {\it uniquely} find a coset 
$s_1(S_u\cap{\rm SU}(3)^\C)$ and another one 
$s_2(\overline{S}_u\cap{\rm SU}(3)^\C)$ such that 
\[g{\rm SU}(3)^\C=\lim\limits_n\big(s_1^{1/n}
(S_u\cap{\rm SU}(3)^\C)^{1/n}\:
s_2^{1/n}(\overline{S}_u\cap{\rm SU}(3)^\C)^{1/n}\big)^n\:\:.\] 
In this way the assignment 
\[t\longmapsto\lim\limits_n\big(\sigma_1(\begin{smallmatrix}
\frac{1}{n}\end{smallmatrix})\sigma_2\begin{smallmatrix}
(\frac{1-t}{n}\end{smallmatrix})\big)^n
\:\Lambda^2\:\lim\limits_n\big(\sigma_1(\begin{smallmatrix}
\frac{1}{n}\end{smallmatrix})\sigma_2(\begin{smallmatrix}\frac{1-t}{n}
\end{smallmatrix})\big)^{-n}\:\:
\lim\limits_n\big(\sigma_1(\begin{smallmatrix}
\frac{1-t}{n}\end{smallmatrix})\sigma_2(\begin{smallmatrix}
\frac{1}{n}\end{smallmatrix})\big)^n
\:\Lambda^2\:\lim\limits_n\big(\sigma_1(\begin{smallmatrix}\frac{1-t}{n}
\end{smallmatrix})\sigma_2(\begin{smallmatrix}\frac{1}{n}
\end{smallmatrix})\big)^{-n}\]
is single-valued and ${\rm SU}(3)^\C$-invariant consequently gives rise to 
a well-defined isotopy connecting 
$g\Lambda^2g^{-1}g\Lambda^2g^{-1}=g\Lambda g^{-1}\in O(\Lambda )^\C$ at $t=0$ 
with $s_1\Lambda^2s_1^{-1}\:s_2\Lambda^2s_2^{-1}\in O'(\Lambda )^\C$ at $t=1$ 
and {\it vice versa}. Clearly $\Lambda\in O(\Lambda )^\C\cap O'(\Lambda )^\C$ 
is kept fixed during this deformation. 

Concerning the assertion itself now, first observe that 
$g=\lim\limits_n\big(s^{1/n}\:\overline{s}^{1/n}\big)^n$ with any $s\in 
S_u$ is a real element i.e. $g\in{\rm G}_2\subset{\rm G}_2^\C$. Thus the 
assertion follows since the homeomorphic image of $O(\Lambda 
)\subset{\rm G}_2$ as in (\ref{palya}) by the isotopy is precisely 
$O'(\Lambda )\subset {\rm G}_2^\C$ given by (\ref{palya'}); however 
$O(\Lambda )$ is homeomorphic to $S^6$ hence so is $O'(\Lambda )$. As a 
global aspect of the construction note that the embedding of $S^6$ 
provided by $O(\Lambda )\subset{\rm G}_2$ hence by 
$O'(\Lambda)\subset{\rm G}_2^\C$ is not null-homotopic (cf. Theorem 
\ref{homotopia}) and obviously the square of the corresponding homotopy 
class is represented by the twin conjugate orbit $O(\Lambda^2)\subset 
{\rm G}_2$ while its third power by $O(\Lambda^3)=e\in{\rm G}_2$ i.e. it 
is already trivial.

After these preliminary constructions and considerations we 
are ready to state

\begin{theorem} 
Take the family ${\mathfrak s}_u\subset{\mathfrak g}_2^\C$ of Samelson 
subalgebras with $u\in P({\mathfrak h}^\C )\setminus P({\mathfrak h})$ and 
consider the induced compact complex $7$-manifolds $Y_u$ which are all 
diffeomorphic to ${\rm G}_2\subset{\rm G}_2^\C$.

Then for every moduli parameter $u$ the corresponding Samelson construction 
gives rise to a compact complex $3$-manifold $X_u$ too which is diffeomorphic 
to $O'(\Lambda )\subset{\rm G}^\C_2$ hence rendering $S^6$ a compact complex 
$3$-manifold.
\label{letezes}
\end{theorem}

\begin{proof}
The strategy is simple: we check that the left-translated and intersected 
Samelson splittings $\big(L_{g*}{\mathfrak s}_u\oplus L_{g*}
\overline{\mathfrak s}_u\big)\cap T_gO'(\Lambda)^\C$ are nice for every 
$g\in O'(\Lambda)$ hence these splittings induce complex structures on 
$T_gO'(\Lambda)$ in the standard way which are integrable along $O'(\Lambda)$.

The standard root basis in ${\mathfrak g}_2^\C$ looks 
like\footnote{An explicit matrix representation of the members of this basis 
will be exhibited soon in Sect. \ref{three}; therefore the interested reader 
can check all assertions about this basis by hand.}
\begin{equation}
\{ H_{\pm a,b}, V_{\pm 1},V_{\pm 2}, V_{\pm 3},U_{\pm 1},
U_{\pm 2},U_{\pm 3}\}
\label{gyok}
\end{equation}
where $a,b$ are real parameters with $a\not=0$. It has the following pleasant 
three properties: the first is that $\{H_{\pm a,b}\}$ span the 
complexified Cartan subalgebra 
${\mathfrak h}^\C\subset{\mathfrak g}_2^\C$; secondly $\{H_{\pm a,b},
V_{\pm k}\}$ with $k=1,2,3$ span 
${\mathfrak s}{\mathfrak u}(3)^\C\subset{\mathfrak g}_2^\C$; 
thirdly $\{ H_{+a,b}, V_{+1}, V_{-2}, V_{-3}, U_{+1}, U_{-2}, U_{+3}\}$ span 
all the family of Samelson subalgebras 
${\mathfrak s}_u\subset{\mathfrak g}_2^\C$ provided 
$P({\mathfrak h}^\C )\setminus P({\mathfrak h})\cong\C P^1\setminus\R P^1
\cong\C\setminus\sqrt{-1}\:\R\ni u =a+\sqrt{-1}\:b$ i.e. $a\not=0$. 
Moreover $\overline{H}_{+a,b}=-H_{-a,b}$,
$\overline{V}_{+k}=V_{-k}$ and $\overline{U}_{+k}=U_{-k}$ hence the
remaining basis elements span the complex conjugate subalgebra
$\overline{\mathfrak s}_u\subset {\mathfrak g}_2^\C$. Finally we 
note that there exist precisely two orthogonal (with respect to the 
$\Ad$-invariant metric on ${\rm G}_2$) complex structures at 
$u_\pm=\pm\frac{\sqrt{3}}{2}-\frac{\sqrt{-1}}{2}$. 

With respect to this basis consider the usual vector space decomposition 
${\mathfrak g}_2^\C={\mathfrak s}{\mathfrak u}(3)^\C\oplus{\mathfrak m}$. 
Intersecting it with ${\mathfrak g}_2^\C={\mathfrak s}_u\oplus 
\overline{\mathfrak s}_u$ we obtain 
${\mathfrak g}_2^\C ={\mathfrak s}_u\cap{\mathfrak s}
{\mathfrak u}(3)^\C\oplus{\mathfrak s}_u\cap
{\mathfrak m}\oplus\overline{\mathfrak s}_u\cap{\mathfrak m}\oplus
{\overline{\mathfrak s}}_u\cap{\mathfrak s}{\mathfrak u}(3)^\C$ 
implying a splitting ${\mathfrak m}=
{\mathfrak s}_u\cap{\mathfrak m}\oplus\overline{\mathfrak s}_u\cap
{\mathfrak m}$. It is straightforward that ${\mathfrak m}$ is spanned by 
$\{U_{\pm k}\}$ with $k=1,2,3$ therefore 
${\mathfrak s}_u\cap{\mathfrak m}$ is spanned by
$\{U_{+1},U_{-2},U_{+3}\}$ and likewise $\overline{\mathfrak s}_u\cap
{\mathfrak m}$ by $\{U_{-1},U_{+2},U_{-3}\}$. Thus we find that 
$\dim_\C\big({\mathfrak s}_u\cap{\mathfrak m}\big)=3=
\dim_\C\big(\overline{\mathfrak s}_u\cap{\mathfrak m}\big)$  
for all moduli parameters $u$.

Let $s_1(t)$ be a smooth curve in 
$S_u$ with the property $s_1(0)=s_1$ and likewise $s_2(t)$ in $\overline{S}_u$ 
such that $s_2(0)=s_2$. Then a generic curve in 
$S_u\Lambda^2S_u^{-1}\subset{\rm G}_2^\C$ passing through 
$s_1\Lambda^2s_1^{-1}$ looks like $s_1(t)\Lambda^2s_1^{-1}(t)$ and similarly 
for the curve $s_2(t)\Lambda^2s_2^{-1}(t)$. Also let 
$L\in{\mathfrak g}_2^\C$ be an element satisfying $\Lambda=\exp L$ hence 
$\Lambda^2=\Lambda^{-1}=\exp(-L)$. Then 
$s_i(t)\Lambda^2s_i^{-1}(t)=\exp(-s_i(t)Ls^{-1}_i(t))$ for $i=1,2$. 
Let us compute the tangent vectors of these curves in the standard way:
\begin{eqnarray}
L^{-1}_{s_i\Lambda^2s_i^{-1}*}\Big(\frac{\dd}{\dd t}\left.
\exp\big(-s_i(t)Ls^{-1}_i(t)\big)\right\vert_{t=0}\Big)&=&
\frac{\dd}{\dd t}\left.\Big(\exp(s_iLs_i^{-1})\exp\big(-s_i(t)Ls^{-1}_i(t)
\big)\Big)\right\vert_{t=0}\nonumber\\
&=&\exp(s_iLs_i^{-1})\exp_*(-s_iLs_i^{-1})\frac{\dd}{\dd t}
\left.\big(-s_i(t)Ls^{-1}_i(t)\big)\right\vert_{t=0}\nonumber\\
&=&\frac{1-{\rm e}^{-\ad_{-s_iLs_i^{-1}}}}{\ad_{-s_iLs_i^{-1}}}
\Big[R_{s_i*}^{-1}\dot{s}_i\:,\:-s_iLs_i^{-1}\Big]\nonumber\\
&=&s_i\left(\frac{{\rm e}^{\ad_L}-1}{\ad_L}\ad_L\big(
L_{s_i*}^{-1}\dot{s}_i\big)\right)s_i^{-1}\nonumber\\
&=&s_i\left(\big({\rm e}^{\ad_L}-1\big)
L_{s_i*}^{-1}\dot{s}_i\right)s_i^{-1}\:\:.\nonumber
\end{eqnarray}
Obviously $L_{s_1*}^{-1}\dot{s}_1\in{\mathfrak s}_u$ and 
$L_{s_2*}^{-1}\dot{s}_2\in\overline{\mathfrak s}_u$. 
Moreover $L\in{\mathfrak h}^\C\subset{\mathfrak g}_2^\C$ 
i.e. it is from the Cartan subalgebra hence $\ad_L$ acts diagonally on 
(\ref{gyok}) thus $\big({\rm e}^{\ad_L}-1\big)L_{s_1*}^{-1}\dot{s}_1
\in{\mathfrak s}_u$ and $\big({\rm e}^{\ad_L}-1\big)L_{s_2*}^{-1}\dot{s}_2
\in\overline{\mathfrak s}_u$ too. The additional datum that 
$\exp L\in Z({\rm SU}(3)^\C)$ therefore 
$\ad_L({\mathfrak s}_u\cap{\mathfrak s}{\mathfrak u}(3)^\C)=0=
\ad_L(\overline{\mathfrak s}_u\cap{\mathfrak s}{\mathfrak u}(3)^\C)$ 
implies that the right hand side in fact belongs to $s_1\big({\mathfrak s}_u
\cap{\mathfrak m}\big)s_1^{-1}\subset{\mathfrak s}_u$ or 
$s_2\big(\overline{\mathfrak s}_u\cap{\mathfrak m}\big)s_2^{-1}\subset
\overline{\mathfrak s}_u$ respectively. 
For simplicity write $p_i:=s_i\Lambda^2s_i^{-1}$ moreover 
$P_1:=S_u\Lambda^2S_u^{-1}$ and $P_2:=
\overline{S}_u\Lambda^2\overline{S}_u^{-1}$. Thus we find that 
$\mbox{$T_{p_1}P_1=L_{p_1*}
\Ad_{s_1}\big({\mathfrak s}_u\cap{\mathfrak m}\big)$\:\:and\:\:
$T_{p_2}P_2=L_{p_2*}\Ad_{s_2}\big(\overline{\mathfrak s}_u
\cap{\mathfrak m}\big)$}$. Recalling that $O'(\Lambda)^\C=P_1P_2$ such that 
$P_1\cap P_2$ is discrete we can apply Leibniz's rule and then insert these 
identities to split the tangent space  
\[T_{p_1p_2}O'(\Lambda)^\C=R_{p_2*}T_{p_1}P_1\oplus 
L_{p_1*}T_{p_2}P_2=L_{p_1*}R_{p_2*}\Ad_{s_1}({\mathfrak s}_u\cap
{\mathfrak m})\oplus
L_{p_1*}L_{p_2*}\Ad_{s_2}(\overline{\mathfrak s}_u\cap{\mathfrak m})\:\:.\]
Put $s_1=s=\overline{s}_2$ hence $p_1=p=\overline{p}_2$ yielding 
$T_{p\overline{p}}O'(\Lambda)^\C=L_{p*}R_{\overline{p}*}
\Ad_{s}({\mathfrak s}_u\cap{\mathfrak m})\oplus L_{p\overline{p}*}
\Ad_{\overline{s}}(\overline{\mathfrak s}_u\cap{\mathfrak m})$. From 
(\ref{palya'}) we know that $p\overline{p}\in O'(\Lambda)$ is a real point 
such that the realification of e.g. the complex summand 
$L_{p\overline{p}*}\Ad_{\overline{s}}(\overline{\mathfrak s}_u
\cap{\mathfrak m})\subset T_{p\overline{p}}O'(\Lambda)^\C$ 
is the real tangent space $T_{p\overline{p}}O'(\Lambda)$ whose 
complexification gives back again $T_{p\overline{p}}O'(\Lambda)^\C$. 
For a complex subspace $V\subseteqq T_g{\rm G}_2^\C$ writing 
its complex conjugate within $T_g{\rm G}_2^\C$ as 
$\overline{V}^g:=L_{g*}\overline{L_{g*}^{-1}V}$ then 
with $V:=L_{p\overline{p}*}\Ad_{\overline{s}}(\overline{\mathfrak s}_u
\cap{\mathfrak m})$ we obtain $T_{p\overline{p}}O'(\Lambda)^\C=(V^\R)^\C=
\overline{V}^{p\overline{p}}\oplus V$. Therefore the asymmetric splitting 
constructed above in the real case induces a new and symmetric one 
\[T_{p\overline{p}}O'(\Lambda)^\C=
L_{p\overline{p}*}\overline{L^{-1}_{p\overline{p}*}L_{p\overline{p}*}
\Ad_{\overline{s}}(\overline{\mathfrak s}_u\cap{\mathfrak m})}
\oplus L_{p\overline{p}*}\Ad_{\overline{s}}
(\overline{\mathfrak s}_u\cap{\mathfrak m})=
L_{p\overline{p}*}\Ad_{s}({\mathfrak s}_u\cap{\mathfrak m})\oplus
L_{p\overline{p}*}\Ad_{\overline{s}}
(\overline{\mathfrak s}_u\cap{\mathfrak m})\] 
moreover readily 
$L_{p\overline{p}*}\Ad_s({\mathfrak s}_u\cap{\mathfrak m})=
L_{p\overline{p}*}{\mathfrak s}_u\cap T_{p\overline{p}}O'(\Lambda)^\C$ and 
$L_{p\overline{p}*}\Ad_{\overline{s}}(\overline{\mathfrak s}_u
\cap{\mathfrak m})=L_{p\overline{p}*}\overline{\mathfrak s}_u\cap 
T_{p\overline{p}}O'(\Lambda)^\C$. 
Thus summing up all of our findings so far we conclude that 
\begin{equation}
T_gO'(\Lambda)^\C=L_{g*}{\mathfrak s}_u\cap T_gO'(\Lambda)^\C\oplus
L_{g*}\overline{\mathfrak s}_u\cap T_gO'(\Lambda)^\C
\label{hasitas}
\end{equation}
is a left-invariant decomposition with 
$\dim_\C (L_{g*}{\mathfrak s}_u\cap T_gO'(\Lambda)^\C )=3=
\dim_\C (L_{g*}\overline{\mathfrak s}_u\cap T_gO'(\Lambda)^\C )$ 
over all $g\in O'(\Lambda )$ and for all Samelson moduli parameters $u
\in P({\mathfrak h}^\C)\setminus P({\mathfrak h})$. 

Thus working over $g\in O'(\Lambda)$ as defined in (\ref{palya'}) 
the Samelson splittings 
$T_g{\rm G}_2^\C =L_{g*}{\mathfrak s}_u\oplus L_{g*}\overline{\mathfrak s}_u$ 
induce sub-splittings $T_gO'(\Lambda )^\C=
Z_{u,g}\oplus\overline{Z}_{u,g}$ where $Z_{u,g}:=L_{g*}{\mathfrak s}_u\cap 
T_gO'(\Lambda )^\C$ and $\overline{Z}_{u,g}:=L_{g*}\overline{\mathfrak s}_u
\cap T_gO'(\Lambda )^\C$. We know additionally that they are both $3$ 
complex dimensional moreover 
$T_gO'(\Lambda)^\C\cong T_gO'(\Lambda)\otimes_\R\C$ and 
$T_gO'(\Lambda)\cap Z_{u,g}=0$. 
Consequently $T_gO'(\Lambda)^\C=Z_{u,g}\oplus\overline{Z}_{u,g}$ gives rise to 
a complex vector space structure $J_{u,g}:T_gO'(\Lambda)\rightarrow 
T_gO'(\Lambda)$ on the underlying real vector space by the general theory. In 
this way $O'(\Lambda)$ is improved to an almost complex manifold 
$(O'(\Lambda ),J_u)$. Concerning its integrability, let $X$ be a real vector 
field along $O'(\Lambda )$ and $X^{1,0}=\frac{1}{2}(X-\sqrt{-1}J_uX)$ its 
corresponding $(1,0)$-type vector field along $(O'(\Lambda ),J_u)$. 
Moreover let $\{Z^{1,0}_1,\dots,Z^{1,0}_7\}$ be a basis over $\C$ of 
${\mathfrak g}_2^{1,0}\subset {\mathfrak g}_2^\C$ hence 
satisfying $\big[Z^{1,0}_i,Z^{1,0}_j\big]^{0,1}=0$ for all $i,j=1,\dots,7$ 
(recall that ${\mathfrak g}_2^\C={\mathfrak g}_2^{1,0}\oplus
{\mathfrak g}_2^{0,1}$ such that ${\mathfrak g}_2^{1,0}={\mathfrak s}_u$ and 
${\mathfrak g}_2^{0,1}=\overline{{\mathfrak s}}_u$). After identifying these 
basis elements with pointwise $\C$-linearly independent left-invariant 
$(1,0)$-type complex vector fields along ${\rm G}_2^\C$ we can pick 
$\C$-valued smooth functions $f_1,\dots, f_7$ along $O'(\Lambda )$ 
(and extend them by zero over the whole ${\rm G}_2^\C$) to write 
$X^{1,0}=\sum_if_iZ^{1,0}_i$. The Nijenhuis bracket of any 
two real vector fields $X,Y\in C^\infty (O'(\Lambda );TO'(\Lambda ))$ then 
looks like 
\[N_{J_u}(X,Y)=\left[X^{1,0}\:,\:Y^{1,0}\right]^{0,1}=\left[\:\:
\sum\limits_{i=1}^7f_iZ^{1,0}_i\:,\:
\sum\limits_{j=1}^7g_jZ^{1,0}_j\right]^{0,1}=\sum\limits_{i,j=1}^7
f_ig_j\left[Z^{1,0}_i\:,\:Z^{1,0}_j\right]^{0,1}=0\:\:.\]  
Thus the Nijenhuis tensor $N_{J_u}$ of $(O'(\Lambda ),J_u)$ itself 
vanishes consequently the almost complex structure is integrable yielding 
a complex manifold $X_u$. Finally we observe that $X_u$ is homeomorphic to 
the $6$-sphere. 
\end{proof}

\begin{remark}\rm 1. Note that {\it a priori} the complex manifolds 
$X_u$ might be non-isomorphic for different values of the 
Samelson moduli parameter $u$. However we will see by the aid of the 
explicit construction at the end of Sect. \ref{three} that these complex 
structures do not depend on $u$. Likewise, repeating everything so far 
with the apparently different twin conjugate orbit 
$O(\Lambda^2)\subset{\rm G}_2$ we do not obtain new complex manifolds 
because we will see in Lemma \ref{egyenloseg} that 
$O(\Lambda^2)=O(\Lambda )$ as subsets of ${\rm G}_2$. 

2. On the one hand Gray \cite{gra} found in 1997 that if $X$ was a 
hypothetical complex manifold diffeomorphic to $S^6$ then 
$H^{0,1}(X)\cong\C$ and raised the question how to interpret the 
generator of this cohomology group. In addition Ugarte \cite[Corollary 
3.3]{uga} proved in 2000 that either (i) $H^{1,1}(X)\not\cong 0$, or 
(ii) $H^{1,1}(X)\cong 0$ and $H^{0,2}(X)\not\cong 0$. Therefore 
$H^{0,1}(X_u)\cong\C$ and either $H^{1,1}(X_u)\not\cong 0$, or 
$H^{1,1}(X_u)\cong 0$ and $H^{0,2}(X_u)\not\cong 0$. On the other hand 
Pittie \cite[Proposition 4.5]{pit} calculated in 1988 the complete 
Dolbeault cohomology ring of $Y_u$ and in particular demonstrated that 
$H^{0,1}(Y_u)\cong\C$, $H^{1,1}(Y_u)\cong\C$ and $H^{0,2}(Y_u)\cong 0$ 
for all moduli parameters $u\in\C P^1\setminus\R P^1$. That is, observe that 
there is a resemblance between the lower degree cohomologies of $X_u$ and $Y_u$.
\end{remark}

\noindent The six-sphere as a complex manifold is not homogeneous. 
Consequently blowing it up once in points belonging to different orbits of 
its automorphism group brings to life spaces which are all diffeomorphic to the 
complex projective three-space but not complex-analytically isomorphic to 
each other \cite{huc-keb-pet}. LeBrun calls in \cite{leb} the 
existence of such {\it exotic $\C P^3$'s} a ``minor disaster''. Here we report 
on a further disaster namely the existence of {\it large exotic $\C^3$'s} in a 
similar sense:
 
\begin{lemma}
Let $x_0\in X_u$ be a point and consider the punctured complex 
manifold $X^\times_u:=X_u\setminus\{x_0\}$. Then the space 
$X^\times_u$ is diffeomorphic to $\C^3$ but is not complex-analytically 
isomorphic to it.
\label{egzotikus}
\end{lemma}

\begin{proof}
Obviously $X^\times_u$ is diffeomorphic to 
$S^6\setminus\{x_0\}$ i.e. to $\R^6$ like $\C^3$ does. Let 
$f: X^\times_u\rightarrow\C$ be a holomorphic function. By Hartogs' theorem 
it extends to a holomorphic function $F:X_u\rightarrow\C$. However $F$ must 
be constant \cite{cam-dem-pet} consequently $f$ is constant on $X^\times_u$ 
as well. Since there exists an abundance of non-trivial holomorphic functions 
on $\C^3$ we conclude that $X^\times_u$ and $\C^3$ are not 
complex-analytically isomorphic. 
\end{proof}


\section{Explicit construction}
\label{three}


In this section we calculate the integrable almost complex tensors $J_u$ 
on $O'(\Lambda)$ underlying the complex manifolds $X_u$ 
of Theorem \ref{letezes}. This computation allows one to read off at least 
that they do not depend on the Samelson parameter $u$ hence the complex 
structures of Theorem \ref{letezes} are in fact equivalent. 

The six-sphere as a complex manifold is not easy to grasp. Since $X_u$ 
is a compact space it cannot be embedded into the affine complex space 
$\C^m$ of any dimension; likewise $H^2(X_u;\C )\cong0$ shows that it is 
not K\"ahler--Hodge consequently it does not admit an embedding into any 
projective complex space $\C P^n$. Therefore its realization as a 
complex submanifold of some well-known complex manifold fails. 
Another odd feature is that the algebraic dimension of 
$X_u$ is zero \cite{cam-dem-pet} which means that all global meromorphic 
functions are constant on it consequently the powerful methods of 
complex analysis also fail to say anything here. What we nevertheless 
can try is to seek the integrable almost complex manifold $(O'(\Lambda), 
J_u)$ underlying $X_u$. This indeed works but the result is so 
complicated that we decided not to display it fully here. In spite of 
this we write down carefully the main steps hence the curious (and 
computer-aided) reader can easily reproduce the calculations and face 
their quantitative complexity directly. Finally, before sinking in the 
heavy details, we raise the question whether or not there exists a 
better general way to exhibit $X_u$ in a more compact or comprehensible 
form.

We begin with an explicit construction of the Samelson family $J_u$ with 
$u\in\C\setminus\sqrt{-1}\:\R$ of all integrable almost complex tensor fields 
on ${\rm G}_2$ as we promised in a footnote of Sect. \ref{two}. We also 
promised in another footnote to write down the root basis (\ref{gyok}) 
of the $14$ dimensional ${\mathfrak g}_2^\C$ explicitly. So let us start with 
this. Our representation of the basis is the smallest possible one and is 
provided by the embedding 
${\mathfrak g}_2^\C\subset{\mathfrak s}{\mathfrak o}(7)^\C$ 
therefore is in terms of $7\times 7$ complex skew symmetric matrices. The 
corresponding matrices are orthonormal with respect to the 
Hermitian (thus {\it not} $\Ad$-invariant!) scalar product 
$\langle V,W\rangle^\C:=\tr (V\overline{W}^T)$ on ${\mathfrak g}_2^\C$ and 
look as follows:
\[H_{\pm a,b}=\frac{1}{2\sqrt{a^2+b^2+b+1}}\left(\begin{smallmatrix}
                0 & 0 & 0 & 0 & 0 & 0 & 0 \\
                0 & 0 &\mp\sqrt{-1} & 0 & 0 & 0 & 0 \\
                0 &\pm\sqrt{-1}  & 0 & 0 & 0 & 0 & 0  \\  
                0 & 0 & 0 & 0 &-a\mp\sqrt{-1}\:b & 0 & 0\\
                0 & 0 & 0 &a\pm\sqrt{-1}\:b & 0 & 0 & 0 \\
                0 & 0 & 0 & 0 & 0 & 0 & -a\mp\sqrt{-1}(b+1) \\
                0 & 0 & 0 & 0 & 0 &a\pm\sqrt{-1}(b+1) & 0 
             \end{smallmatrix}\right) 
\:\:\:\begin{array}{ll}\mbox{$a,b\in\R$,}\\
                        \mbox{$a\not=0$}
       \end{array}\]
and one checks that $\{H_{\pm a,b}\}$ span the $2$ dimensional complex Cartan 
subalgebra of ${\mathfrak h}^\C\subset{\mathfrak g}_2^\C$; moreover 
\[V_{\pm 1}=\frac{1}{2\sqrt{2}}\left(\begin{smallmatrix}
                0 & 0 & 0 & 0 & 0 & 0 & 0 \\
                0 & 0 & 0 & -1 & \mp\sqrt{-1} & 0 & 0 \\
                0 & 0 & 0 & \pm\sqrt{-1} & -1 & 0 & 0 \\
                0 & 1 & \mp\sqrt{-1} & 0 & 0 & 0 & 0 \\
                0 & \pm\sqrt{-1} & 1 & 0 & 0 & 0 & 0 \\
                0 & 0 & 0 & 0 & 0 & 0 & 0 \\
                0 & 0 & 0 & 0 & 0 & 0 & 0                  
             \end{smallmatrix}\right)\] 
\[V_{\pm 2}=\frac{1}{2\sqrt{2}}\left(\begin{smallmatrix}
                0 & 0 & 0 & 0 & 0 & 0 & 0 \\
                0 & 0 & 0 & 0 & 0 & -1 &\pm\sqrt{-1} \\
                0 & 0 & 0 & 0 & 0 & \pm\sqrt{-1} & 1 \\
                0 & 0 & 0 & 0 & 0 & 0 & 0 \\
                0 & 0 & 0 & 0 & 0 & 0 & 0 \\
                0 & 1 & \mp\sqrt{-1} & 0 & 0 & 0 & 0 \\
                0 & \mp\sqrt{-1} & -1 & 0 & 0 & 0 &0 
             \end{smallmatrix}\right)\]
\[V_{\pm 3}=\frac{1}{2\sqrt{2}}\left(\begin{smallmatrix}
                 0 & 0 & 0 & 0 & 0 & 0 & 0 \\
                 0 & 0 & 0 & 0 & 0 & 0 & 0 \\
                 0 & 0 & 0 & 0 & 0 & 0 & 0 \\
                 0 & 0 & 0 & 0 & 0 & -1 &\pm\sqrt{-1} \\
                 0 & 0 & 0 & 0 & 0 & \pm\sqrt{-1} & 1 \\
                 0 & 0 & 0 & 1 & \mp\sqrt{-1} & 0 & 0 \\
                 0 & 0 & 0 & \mp\sqrt{-1} & -1 & 0 & 0 
              \end{smallmatrix}\right)\]
and one checks that $\{H_{\pm a,b},V_{\pm 1}, V_{\pm 2}, V_{\pm 3}\}$ span 
the (maximal) subalgebra 
${\mathfrak s}{\mathfrak u}(3)^\C\subset{\mathfrak g}_2^\C$; and finally 
\[U_{\pm 1}=\frac{1}{2\sqrt{6}}\left(\begin{smallmatrix}
              0 &\mp2\sqrt{-1} & -2 & 0 & 0 & 0 & 0 \\
              \pm2\sqrt{-1}& 0 & 0 & 0 & 0 & 0 & 0 \\
              2 & 0 & 0 & 0 & 0 & 0 & 0 \\
              0 & 0 & 0 & 0 & 0 & -1 & \pm\sqrt{-1} \\
              0 & 0 & 0 & 0 & 0 & \mp\sqrt{-1} & -1 \\
              0 & 0 & 0 & 1 & \pm\sqrt{-1} & 0 & 0 \\
              0 & 0 & 0 & \mp\sqrt{-1} & 1 & 0 & 0 \\
           \end{smallmatrix}\right)\]
\[U_{\pm 2}=\frac{1}{2\sqrt{6}}\left(\begin{smallmatrix}
              0 & 0 & 0 & \pm 2\sqrt{-1} & 2 & 0 & 0 \\
              0 & 0 & 0 & 0 & 0 & -1 & \pm\sqrt{-1} \\
              0 & 0 & 0 & 0 & 0 & \mp\sqrt{-1} & -1 \\
              \mp 2\sqrt{-1} & 0 & 0 & 0 & 0 & 0 & 0 \\
              -2 & 0 & 0 & 0 & 0 & 0 & 0 \\
              0 & 1 & \pm\sqrt{-1} & 0 & 0 & 0 & 0 \\
              0 & \mp\sqrt{-1} & 1 & 0 & 0 & 0 & 0
            \end{smallmatrix}\right)\]
\[U_{\pm 3}=\frac{1}{2\sqrt{6}}\left(\begin{smallmatrix}
             0 & 0 & 0 & 0 & 0 & 2 & \mp2\sqrt{-1} \\
             0 & 0 & 0 & \pm\sqrt{-1} & 1 & 0 & 0 \\
             0 & 0 & 0 & 1 & \mp\sqrt{-1} & 0 & 0 \\
             0 & \mp\sqrt{-1} & -1 & 0 & 0 & 0 & 0 \\
             0 & -1 & \pm\sqrt{-1}& 0 & 0 & 0 & 0 \\
             -2 & 0 & 0 & 0 & 0 & 0 & 0 \\
             \pm2\sqrt{-1} & 0 & 0 & 0 & 0 & 0 & 0
           \end{smallmatrix}\right)\]
and one checks that $\{H_{+a,b},V_{+1},V_{-2},V_{-3},U_{+1},U_{-2},U_{+3}\}$ 
span ${\mathfrak s}_u
\subset{\mathfrak g}_2^\C$ where $u\in\C\setminus\sqrt{-1}\:\R$ 
is written in the form $u=a+\sqrt{-1}\:b$ with $a\not=0$; we also find that 
the remaining elements from the root basis 
$\{H_{-a,b}, V_{-1},V_{+2},V_{+3},U_{-1},U_{+2},U_{-3}\}$ form a 
basis in the complex conjugate subalgebra 
$\overline{\mathfrak s}_u\subset{\mathfrak g}_2^\C$. 

Consequently picking $a_i+\sqrt{-1}\:b_i\in\C$ ($i=0,1,\dots,6$) to write 
an element $W\in{\mathfrak s}_u$ as 
\begin{eqnarray}
W&:=&2\sqrt{a^2+b^2+b+1}\:(a_0+\sqrt{-1}\:b_0)H_{+a,b}\nonumber\\
&&+2\sqrt{2}(a_4+\sqrt{-1}\:b_4)V_{+1}+
2\sqrt{2}(a_5+\sqrt{-1}\:b_5)V_{-2}+2\sqrt{2}
(a_6+\sqrt{-1}\:b_6)V_{-3}\nonumber\\
&&+2\sqrt{6}(a_1+\sqrt{-1}\:b_1)U_{+1}+
2\sqrt{6}(a_2+\sqrt{-1}\:b_2)U_{-2}+2\sqrt{6}(a_3+\sqrt{-1}\:b_3)U_{+3}\nonumber
\end{eqnarray}
and putting $J_u:{\mathfrak g}_2\rightarrow{\mathfrak g}_2$ to be 
$J_u(\re W):=-\im W$ dictated by the general theory we 
obtain an $\R$-linear transformation $J_u\in{\rm End}\:{\mathfrak g}_2$. 
Its action on 
\[\re W=\left(\begin{smallmatrix}
0&2b_1&-2a_1&-2b_2&2a_2&2a_3&2b_3\\
-2b_1&0&\frac{1}{2}b_0&-b_3-a_4&a_3+b_4&-a_2-a_5&-b_2-b_5\\
2a_1&-\frac{1}{2}b_0&0&a_3-b_4&-a_4+b_3&b_2-b_5&-a_2+a_5\\
2b_2&a_4+b_3&-a_3+b_4&0&-\frac{1}{2}aa_0+\frac{1}{2}bb_0&-a_1-a_6&-b_1-b_6\\
-2a_2&-a_3-b_4&a_4-b_3&\frac{1}{2}aa_0-\frac{1}{2}bb_0&0&b_1-b_6&-a_1+a_6\\
-2a_3&a_2+a_5&-b_2+b_5&a_1+a_6&-b_1+b_6&0&-\frac{1}{2}aa_0+\frac{1}{2}(b+1)b_0\\
-2b_3&b_2+b_5&a_2-a_5&b_1+b_6&a_1-a_6&\frac{1}{2}aa_0-\frac{1}{2}(b+1)b_0&0
\end{smallmatrix}\right)\]
is by definition 
\[-\im W=\left(\begin{smallmatrix}
0&2a_1&2b_1&-2a_2&-2b_2&-2b_3&2a_3\\
-2a_1&0&\frac{1}{2}a_0&-a_3+b_4&a_4-b_3&b_2+b_5&-a_2-a_5\\
-2b_1&-\frac{1}{2}a_0&0&-a_4-b_3&a_3+b_4&a_2-a_5&b_2-b_5\\
2a_2&a_3-b_4&a_4+b_3&0&\frac{1}{2}ba_0+\frac{1}{2}ab_0&b_1+b_6&-a_1-a_6\\
2b_2&-a_4+b_3&-a_3-b_4&-\frac{1}{2}ba_0-\frac{1}{2}ab_0&0&a_1-a_6&b_1-b_6\\
2b_3&-b_2-b_5&-a_2+a_5&-b_1-b_6&-a_1+a_6&0&\frac{1}{2}(b+1)a_0+\frac{1}{2}ab_0\\
-2a_3&a_2+a_5&-b_2+b_5&a_1+a_6&-b_1+b_6&-\frac{1}{2}(b+1)a_0-\frac{1}{2}ab_0&0
\end{smallmatrix}\right)\]
hence we immediately check that $J_u^2=-\id_{{\mathfrak g}_2}$. The 
shape of $J_u$ can be read off from these matrices more explicitly if 
we introduce the real orthonormal basis 
\begin{eqnarray}
H_+&:=&\frac{\sqrt{a^2+b^2+b+1}}{2a}(H_{+a,b}+H_{-a,b})\nonumber\\
H_-&:=&\frac{\sqrt{a^2+b^2+b+1}}{\sqrt{-3}}(H_{+a,b}-H_{-a,b})
-\frac{2b+1}{\sqrt{3}}\:H_+\nonumber\\
X_{\pm k}&:=&\frac{1}{\sqrt{\pm 2}}(U_{+k}\pm
U_{-k})\:\:\:,\:\:\:k=1,2,3\nonumber\\
Y_{\pm k}&:=&\frac{1}{\sqrt{\pm 2}}(V_{+k}\pm
V_{-k})\:\:\:,\:\:\:k=1,2,3\nonumber
\end{eqnarray}
(note that $H_\pm$ are already independent of $a,b$) on the compact 
real form ${\mathfrak g}_2\subset{\mathfrak s}{\mathfrak o}(7)$ of 
${\mathfrak g}_2^\C$ equipped with the (positive definite) 
$\Ad$-invariant real scalar product $\langle V,W\rangle:=\tr( V\:W^T)$. 
A straightforward computation verifies that in this basis the action of 
$J_u$ takes the simple blockdiagonal shape 
\begin{eqnarray}
J_uH_+=-\frac{2b+1}{2a}H_++
\frac{4a^2+4b^2+4b+1}{2\sqrt{3}\:a}H_-&,&
J_uH_-=-\frac{\sqrt{3}}{2a}H_++\frac{2b+1}{2a}H_-\nonumber\\
J_uX_{+k}=X_{-k}&,&J_uX_{-k}=-X_{+k}\label{Jhatas}\\
J_uY_{+k}=Y_{-k}&,&J_uY_{-k}=-Y_{+k}\:\:.\nonumber
\end{eqnarray}
One can find precisely two integrable almost complex structures $J_{u_{\pm}}$ 
at $a=\pm\frac{\sqrt{3}}{2}$ and $b=-\frac{1}{2}$ which are orthogonal for the 
aforementioned natural $\Ad$-invariant real scalar product on 
${\mathfrak g}_2$ i.e. $\langle J_{u_\pm} V,V\rangle=0$ and 
$\vert J_{u_\pm}V\vert=\vert V\vert$ for all $V\in{\mathfrak g}_2$. After 
left-translating the $J_u$'s over the whole compact group they give rise to 
complex structures on ${\rm G}_2$ such that the two orthogonal $J_{u_\pm}$ 
yield two complex structures compatible with the bi-invariant 
metric on ${\rm G}_2$. This completes the construction of all the integrable 
almost complex tensors $J_u$ on ${\rm G}_2$ {\it \`a la} Samelson 
(for further details cf. \cite[Example on p. 123]{pit}). 

Our next task is to construct the deformed conjugate orbit $O'(\Lambda)\subset
{\rm G}_2^\C$. This will be achieved in two technical 
steps below. First we construct the original conjugate orbit 
(\ref{palya}) secondly its deformation as defined in (\ref{palya'}). 

Recall the classical fact that ${\rm G}_2$ coincides with the automorphism 
group of the octonions over the reals; this important fact has not been used 
so far explicitly. We shall identify $O(\Lambda )\subset{\rm G}_2$ with the 
subset of inner automorphisms of the octonions. Let $\OO$ denote the 
non-associative, 
unital normed algebra of the octonions (or Cayley numbers) over the reals. In 
the canonical oriented basis $\{\ee_0,\ee_1,\dots,\ee_7\}$ the $\ee_0$ plays 
the role of the unit hence $\re\OO:=\R\ee_0\subset\OO$ is the real part. To be 
absolutely unambiguous we communicate our octonionic multiplication convention 
here:
\[\begin{tabular}{c|cccccccc}
  &$\ee_0$&$\ee_1$&$\ee_2$&$\ee_3$&$\ee_4$&$\ee_5$&$\ee_6$&$\ee_7$ \\
\hline
$\ee_0$&$\ee_0$&$\ee_1$&$\ee_2$&$\ee_3$&$\ee_4$&$\ee_5$&$\ee_6$&$\ee_7$\\
$\ee_1$&$\ee_1$&$-\ee_0$&$\ee_3$&$-\ee_2$&$\ee_5$&$-\ee_4$&$-\ee_7$&$\ee_6$\\
$\ee_2$&$\ee_2$&$-\ee_3$&$-\ee_0$&$\ee_1$&$\ee_6$&$\ee_7$&$-\ee_4$&$-\ee_5$\\
$\ee_3$&$\ee_3$&$\ee_2$&$-\ee_1$&$-\ee_0$&$\ee_7$&$-\ee_6$&$\ee_5$&$-\ee_4$\\
$\ee_4$&$\ee_4$&$-\ee_5$&$-\ee_6$&$-\ee_7$&$-\ee_0$&$\ee_1$&$\ee_2$&$\ee_3$\\
$\ee_5$&$\ee_5$&$\ee_4$&$-\ee_7$&$\ee_6$&$-\ee_1$&$-\ee_0$&$-\ee_3$&$\ee_2$\\
$\ee_6$&$\ee_6$&$\ee_7$&$\ee_4$&$-\ee_5$&$-\ee_2$&$\ee_3$&$-\ee_0$&$-\ee_1$\\
$\ee_7$&$\ee_7$&$-\ee_6$&$\ee_5$&$\ee_4$&$-\ee_3$&$-\ee_2$&$\ee_1$&$-\ee_0$
\end{tabular}\]
(actually there are many different conventions in use). The basis gives 
rise to a canonical $\R$-linear isomorphism of oriented spaces 
$(\OO ,\ee_0,\dots,\ee_7)\cong\R^8$. We can use the standard scalar 
product on $\R^8$ to define $\im\OO:=(\re\OO )^\perp\subset\OO$ and 
introduce a multiplicative norm $\vert\:\cdot\:\vert$ on $\OO$. Then 
canonically $(\im\OO ,\ee_1,\dots,\ee_7)\cong\R^7$ 
and in this way we can look at the six-sphere as the set of imaginary 
octonions of unit length i.e. we will suppose $S^6\subset\im\OO$. If 
${\bf u},{\bf v}\in\OO$ and ${\bf v}\not=0$ then the identity of elasticity 
convinces us that $({\bf v}{\bf u}){\bf v}^{-1}={\bf v}({\bf u}{\bf v}^{-1})$ 
hence it is meaningful to talk about {\it inner 
automorphisms} of the octonions. An important result \cite{lam} says that 
${\bf v}$ indeed gives rise to an inner automorphism if and only if 
$4(\re{\bf v})^2=\vert{\bf v}\vert^2$ i.e. $3(\re {\bf v})^2=\vert\im{\bf v}
\vert^2$ holds. Note that this condition implies ${\bf v}^3$ is a non-zero 
scalar therefore it corresponds to the trivial automorphism of $\OO$. Picking 
an ${\bf x}\in S^6$ the non-real octonion ${\bf v}:=\ee_0+\sqrt{3}\:{\bf x}$ 
satisfying ${\bf v}^3=-8\ee_0$ therefore gives an inner automorphism 
(hence parantheses can be omitted)
\begin{equation}
{\bf u}\longmapsto (\ee_0+\sqrt{3}\:{\bf x}){\bf u}(\ee_0+\sqrt{3}
\:{\bf x})^{-1}=
\frac{1}{4}(\ee_0+\sqrt{3}\:{\bf x}){\bf u}(\ee_0-\sqrt{3}\:{\bf x})
\label{konjugalas}
\end{equation} 
and all inner automorphisms of the octonions are of this form. In this way 
we get a remarkable map 
\begin{equation}
f:S^6\longrightarrow {\rm G}_2\:\:\:.
\label{belso.automorfizmusok}
\end{equation}
Knowing that $\re\OO$ is invariant under all automorphisms and that the 
corresponding reduced linear map of $\im\OO$ is an orientation preserving 
orthogonal transformation of $\R^7$ we can embed ${\rm G}_2$ into 
${\rm SO}(7)$ as usual. Under the canonical isomorphism 
$(\im\OO, \ee_1,\dots,\ee_7)\cong\R^7$ put 
${\bf x}=x_1\ee_1 +\dots +x_7\ee_7$ and thus the image $f({\bf x})$ of the map 
(\ref{belso.automorfizmusok}) at ${\bf x}$ takes the impressive shape 
\[f\left(\begin{smallmatrix}
x_1\\
x_2\\
x_3\\
x_4\\
x_5\\
x_6\\
x_7
\end{smallmatrix}\right)
=\frac{\sqrt{3}}{2}\left(\begin{smallmatrix}
-\frac{1}{\sqrt{3}}+\sqrt{3}\:x_1^2&-x_3+\sqrt{3}\:x_1x_2&x_2+\sqrt{3}\:x_1x_3
&-x_5+\sqrt{3}\:x_1x_4&x_4+\sqrt{3}\:x_1x_5&x_7+\sqrt{3}\:x_1x_6&
-x_6+\sqrt{3}\:x_1x_7\\
x_3+\sqrt{3}\:x_2x_1&-\frac{1}{\sqrt{3}}+\sqrt{3}\:x_2^2&-x_1+\sqrt{3}\:x_2x_3
&-x_6+\sqrt{3}\:x_2x_4&-x_7+\sqrt{3}\:x_2x_5&x_4+\sqrt{3}\:x_2x_6&
x_5+\sqrt{3}\:x_2x_7\\
-x_2+\sqrt{3}\:x_3x_1&x_1+\sqrt{3}\:x_3x_2&-\frac{1}{\sqrt{3}}+\sqrt{3}\:x_3^2&
-x_7+\sqrt{3}\:x_3x_4&x_6+\sqrt{3}\:x_3x_5&-x_5+\sqrt{3}\:x_3x_6
&x_4+\sqrt{3}\:x_3x_7\\
x_5+\sqrt{3}\:x_4x_1&x_6+\sqrt{3}\:x_4x_2&x_7+\sqrt{3}\:x_4x_3
&-\frac{1}{\sqrt{3}}+\sqrt{3}\:x_4^2&-x_1+\sqrt{3}\:x_4x_5&-x_2+
\sqrt{3}\:x_4x_6&-x_3+\sqrt{3}\:x_4x_7\\
-x_4+\sqrt{3}\:x_5x_1&x_7+\sqrt{3}\:x_5x_2&-x_6+\sqrt{3}\:x_5x_3
&x_1+\sqrt{3}\:x_5x_4&-\frac{1}{\sqrt{3}}+\sqrt{3}\:x_5^2&x_3+\sqrt{3}\:x_5x_6&
-x_2+\sqrt{3}\:x_5x_7\\
-x_7+\sqrt{3}\:x_6x_1&-x_4+\sqrt{3}\:x_6x_2&x_5+\sqrt{3}\:x_6x_3
&x_2+\sqrt{3}\:x_6x_4&-x_3+\sqrt{3}\:x_6x_5&-\frac{1}{\sqrt{3}}+
\sqrt{3}\:x_6^2&x_1+\sqrt{3}\:x_6x_7\\
x_6+\sqrt{3}\:x_7x_1&-x_5+\sqrt{3}\:x_7x_2&-x_4+\sqrt{3}\:x_7x_3&
x_3+\sqrt{3}\:x_7x_4&x_2+\sqrt{3}\:x_7x_5&-x_1+\sqrt{3}\:x_7x_6&
-\frac{1}{\sqrt{3}}+\sqrt{3}\:x_7^2
\end{smallmatrix}\right)\]
of a $7\times 7$ special orthogonal matrix with 
$x_1,\dots ,x_7\in\R$ satisfying $x_1^2+\dots +x_7^2=1$. We shall need 
the derivative 
$f({\bf x})_*:T_{\bf x}S^6\rightarrow T_{f({\bf x})}{\rm G}_2$ 
at ${\bf x}\in S^6$ of this map, too. 
If $\xi\in T_{\bf x}S^6$ is a tangent vector then $f({\bf x})_*\xi\in 
T_{f({\bf x})}{\rm G}_2$ is
its image. Putting $\xi=\xi_1\ee_1+\dots +\xi_7\ee_7$ satisfying 
$x_1\xi_1+\dots+x_7\xi_7=0$ a long but 
straightforward computation verifies that
\begin{eqnarray}
f({\bf x})_*\left(\begin{smallmatrix}
\xi_1 \\
\xi_2 \\
\xi_3 \\
\xi_4 \\
\xi_5 \\
\xi_6 \\
\xi_7\end{smallmatrix}\right)
\!\!\!&=&\!\!\!
\frac{3}{2}\left(\begin{smallmatrix}
2x_1\xi_1 & -\frac{1}{\sqrt{3}}\xi_3+x_1\xi_2+x_2\xi_1 &
\frac{1}{\sqrt{3}}\xi_2+x_1\xi_3+x_3\xi_1 &
-\frac{1}{\sqrt{3}}\xi_5+x_1\xi_4+x_4\xi_1 &
\frac{1}{\sqrt{3}}\xi_4+x_1\xi_5+x_5\xi_1 \\
\frac{1}{\sqrt{3}}\xi_3+x_2\xi_1+x_1\xi_2 & 2x_2\xi_2 &
-\frac{1}{\sqrt{3}}\xi_1+x_2\xi_3+x_3\xi_2 &
-\frac{1}{\sqrt{3}}\xi_6+x_2\xi_4+x_4\xi_2 &
-\frac{1}{\sqrt{3}}\xi_7+x_2\xi_5+x_5\xi_2 \\
-\frac{1}{\sqrt{3}}\xi_2+x_3\xi_1+x_1\xi_3 &
\frac{1}{\sqrt{3}}\xi_1+x_3\xi_2+x_2\xi_3 & 2x_3\xi_3 &
-\frac{1}{\sqrt{3}}\xi_7+x_3\xi_4+x_4\xi_3 &
\frac{1}{\sqrt{3}}\xi_6+x_3\xi_5+x_5\xi_3 \\
\frac{1}{\sqrt{3}}\xi_5+x_4\xi_1+x_1\xi_4 &
\frac{1}{\sqrt{3}}\xi_6+x_4\xi_2+x_2\xi_4 &
\frac{1}{\sqrt{3}}\xi_7+x_4\xi_3+x_3\xi_4 & 2x_4\xi_4 &
-\frac{1}{\sqrt{3}}\xi_1+x_4\xi_5+x_5\xi_4 \\
-\frac{1}{\sqrt{3}}\xi_4+x_5\xi_1+x_1\xi_5 &
\frac{1}{\sqrt{3}}\xi_7+x_5\xi_2+x_2\xi_5 &
-\frac{1}{\sqrt{3}}\xi_6+x_5\xi_3+x_3\xi_5 &
\frac{1}{\sqrt{3}}\xi_1+x_5\xi_4+x_4\xi_5 & 2x_5\xi_5 \\
-\frac{1}{\sqrt{3}}\xi_7+x_6\xi_1+x_1\xi_6 &
-\frac{1}{\sqrt{3}}\xi_4+x_6\xi_2+x_2\xi_6 &
\frac{1}{\sqrt{3}}\xi_5+x_6\xi_3+x_3\xi_6 &
\frac{1}{\sqrt{3}}\xi_2+x_6\xi_4+x_4\xi_6 &
-\frac{1}{\sqrt{3}}\xi_3+x_6\xi_5+x_5\xi_6 \\
\frac{1}{\sqrt{3}}\xi_6+x_7\xi_1+x_1\xi_7 &
-\frac{1}{\sqrt{3}}\xi_6+x_7\xi_2+x_2\xi_7 &
-\frac{1}{\sqrt{3}}\xi_4+x_7\xi_3+x_3\xi_7 &
\frac{1}{\sqrt{3}}\xi_3+x_7\xi_4+x_4\xi_7 &
\frac{1}{\sqrt{3}}\xi_2+x_7\xi_5+x_5\xi_7
 \end{smallmatrix}\right.\nonumber\\
& &\nonumber\\
& &\hspace{8.2cm}\left.\begin{smallmatrix}
\frac{1}{\sqrt{3}}\xi_7+x_1\xi_6+x_6\xi_1 &
-\frac{1}{\sqrt{3}}\xi_6+x_1\xi_7+x_7\xi_1\\
\frac{1}{\sqrt{3}}\xi_4+x_2\xi_6+x_6\xi_2 &
\frac{1}{\sqrt{3}}\xi_5+x_2\xi_7+x_7\xi_2 \\
-\frac{1}{\sqrt{3}}\xi_5+x_3\xi_6+x_6\xi_3 &
\frac{1}{\sqrt{3}}\xi_4+x_3\xi_7+x_7\xi_3\\
-\frac{1}{\sqrt{3}}\xi_2+x_4\xi_6+x_6\xi_4 &
-\frac{1}{\sqrt{3}}\xi_3+x_4\xi_7+x_7\xi_4 \\
\frac{1}{\sqrt{3}}\xi_3+x_5\xi_6+x_6\xi_5 &
-\frac{1}{\sqrt{3}}\xi_2+x_5\xi_7+x_7\xi_5 \\
2x_6\xi_6 & \frac{1}{\sqrt{3}}\xi_1+x_6\xi_7+x_7\xi_6 \\
-\frac{1}{\sqrt{3}}\xi_1+x_7\xi_6+x_6\xi_7 & 2x_7\xi_7
\end{smallmatrix}\right)\nonumber
\end{eqnarray}
out of which we also obtain 
(but already will be unable to plot anything from now on) the pullback 
$L_{f({\bf x})\:*}^{-1}(f({\bf x})_*\xi )\in T_e{\rm G}_2=
{\mathfrak g}_2$. In less fancy notation this is $f({\bf x})^{-1}
f({\bf x})_*\xi\in\R (7)$ 
i.e. the plain matrix product of the group-theoretic inverse 
$f({\bf x})^{-1}$ and $f({\bf x})_*\xi$ above.

It readily follows from (\ref{konjugalas}) or by a direct computation that 
$f({\bf x})^3=1_{\R^7}$ and $f({\bf x}){\bf x}={\bf x}$ i.e. we can visualize 
$f({\bf x})$ as a degree $\frac{2\pi}{3}$ rotation $R_{\bf x}$ about the axis 
through ${\bf x}\in\R^7$. Since $\Lambda\in Z({\rm SU}(3))\cong\Z_3$ satisfies 
$\Lambda^3=e\in{\rm G}_2$ it is also true that $h^3=e$ for all 
$h\in O(\Lambda )$. In fact the two subsets $f(S^6)$ and $O(\Lambda )$ of 
${\rm G}_2$ are nothing but the same:
\begin{lemma}{\rm (cf. \cite[pp. 160-161]{cha-rig})}  
The conjugate orbit $O(\Lambda )\subset{\rm G}_2$ of (\ref{palya}) 
and the image $f(S^6)\subset{\rm G}_2$ of the map (\ref{belso.automorfizmusok}) 
coincide as subsets within ${\rm G}_2$ i.e. 
\[O(\Lambda )=\big\{f({\bf x})\:\big\vert\:{\bf x}\in S^6\big\}\:\:.\] 
Moreover $O(\Lambda )=O(\Lambda^2)$, where 
$O(\Lambda^2)=\{g\Lambda^2g^{-1}\vert g\in{\rm G}_2\}$ is the conjugate 
orbit passing through the square of the generator $\Lambda^2\in 
Z({\rm SU}(3))\subset{\rm G}_2$. 
\label{egyenloseg} 
\end{lemma}

\begin{proof} 
We quickly observe that
\begin{equation}
f\left(\begin{smallmatrix}
1\\
0\\
0\\
0\\
0\\
0\\
0
\end{smallmatrix}\right)=\left(\begin{smallmatrix}
1 & 0 & 0 & 0 & 0 & 0 & 0 \\
0 & -\frac{1}{2} & -\frac{\sqrt{3}}{2} & 0 & 0 & 0 & 0\\
0 & \frac{\sqrt{3}}{2} & -\frac{1}{2}  & 0 & 0 & 0 & 0\\
0 & 0 & 0 & -\frac{1}{2} & -\frac{\sqrt{3}}{2} & 0 & 0\\
0 & 0 & 0 & \frac{\sqrt{3}}{2} & -\frac{1}{2}  & 0 & 0\\
0 & 0 & 0 & 0 & 0 & -\frac{1}{2} & \frac{\sqrt{3}}{2} \\
0 & 0 & 0 & 0 & 0 & -\frac{\sqrt{3}}{2} & -\frac{1}{2}
\end{smallmatrix}\right)
\label{Lambda}
\end{equation}
hence $f(\ee_1)=\Lambda\in Z({\rm SU}(3))\subset{\rm G}_2$. 
Therefore the action ${\bf u}\mapsto\Lambda{\bf u}$ on ${\bf u}\in\OO$ 
arises from the inner automorphism ${\bf u}\mapsto (\ee_0+\sqrt{3}\:\ee_1)
{\bf u}(\ee_0+\sqrt{3}\:\ee_1)^{-1}$ 
in (\ref{konjugalas}). Now pick $g\in{\rm G}_2$ then the twisted action 
${\bf u}\mapsto (g\Lambda g^{-1}){\bf u}$ looks like 
\begin{eqnarray}
(g\Lambda g^{-1}){\bf u}\!\!\!&=&\!\!\!g\left((\ee_0+\sqrt{3}\:\ee_1)(g^{-1}
{\bf u})(\ee_0+\sqrt{3}\:\ee_1)^{-1}\right)
=-\frac{1}{8}\:g\!\!\left((\ee_0+\sqrt{3}\:\ee_1)(g^{-1}{\bf u})
(\ee_0+\sqrt{3}\:\ee_1)^2\!\right)\nonumber\\
&=&-\frac{1}{8}\:g(\ee_0+\sqrt{3}\:\ee_1){\bf u}(g(\ee_0+\sqrt{3}\:\ee_1))^2=
-\frac{1}{8}(\ee_0+\sqrt{3}(g\ee_1)){\bf u}(\ee_0+\sqrt{3}(g\ee_1))^2
\nonumber\\
&=&(\ee_0+\sqrt{3}(g\ee_1)){\bf u}(\ee_0+\sqrt{3}(g\ee_1))^{-1}\nonumber
\end{eqnarray}
consequently it comes from an inner automorphism by 
${\bf v}:=\ee_0+\sqrt{3}(g\ee_1)$. Therefore $f(g\ee_1)=g\Lambda g^{-1}$ and 
taking into account that ${\rm G}_2$ acts transitively on 
$S^6\subset\im\OO$ (with stabilizer subgroup 
${\rm SU}(3)\subset{\rm G}_2$) we conclude that $f(S^6)=O(\Lambda )$. 

The fact $f(-{\bf x})=f({\bf x})^T=f({\bf x})^{-1}$ gives the identity 
$f(\pm{\bf x})=f({\bf x})^{\pm 1}$ 
for all ${\bf x}\in S^6$. Therefore $f(\pm\ee_1)=f(\ee_1 )^{\pm 1}=
\Lambda^{\pm 1}$ yielding $f(-\ee_1)=\Lambda^{-1}=\Lambda^2$. Taking an element 
$g\in{\rm G}_2$ satisfying $g\ee_1=-\ee_1$ (unique up to two-sided 
multiplication with elements of ${\rm SU}(3)\subset{\rm G}_2$) we can write 
$\Lambda^2=g\Lambda g^{-1}$; hence the two conjugate orbits of $\Lambda$ 
and $\Lambda^2$ in ${\rm G}_2$ are not distinct consequently they 
must coincide. 
\end{proof}

\noindent Now we are in a position to construct the deformed conjugate orbit 
$O'(\Lambda )\subset{\rm G}_2^\C$. Replace the real vector 
${\bf x}=x_1{\bf e}_1+\dots +x_7{\bf e}_7$ satisfying 
$x_1^2+\dots+x_7^2=1$ with a complex vector 
${\bf z}=z_1{\bf e}_1+\dots +z_7{\bf e}_7$ such that 
$z_1^2+\dots+z_7^2=1$. Then ${\bf w}={\bf e}_0+\sqrt{3}\:{\bf z}$ continues to 
be invertible and satisfies $3(\re{\bf w})^2=\vert\im{\bf w}
\vert^2$ hence generates an inner automorphism of the complexified 
octonions $\OO\otimes_\R\C$. Consequently inserting ${\bf z}$ into 
(\ref{belso.automorfizmusok}) we obtain the complexified map 
\[f^\C:(S^6)^\C\longrightarrow{\rm G}_2^\C\] 
and by Lemma \ref{egyenloseg} we know that $O(\Lambda)^\C=\{f^\C({\bf z})
\:\vert\:{\bf z}\in (S^6)^\C\}$ and $O(\Lambda)^\C=O(\Lambda^2)^\C$. 

\begin{lemma} There exists a bounded domain $\Omega\subset\C^3$ about the 
origin (whose more precise shape is not important) 
such that the deformed conjugate orbit $O'(\Lambda)\subset{\rm G}_2^\C$ as 
defined in (\ref{palya'}) looks like 
\[\begin{array}{lll}
O'(\Lambda)=\Big\{f^\C({\bf z})f^\C(\overline{{\bf z}}) & \Big\vert &
      {\bf z}= -{\bf e}_1+z_2{\bf e}_2-\sqrt{-1}z_2{\bf e}_3+z_4{\bf e}_4
     -\sqrt{-1}z_4{\bf e}_5+z_6{\bf e}_6-\sqrt{-1}z_6{\bf e}_7\:\:,\\
              & & (z_2,z_4,z_6)\in\overline{\Omega}\Big\}
\end{array}\]
where $\overline{\bf z}$ denotes the conjugate of ${\bf z}$ as a 
complex vector in $\C^7$ (and not the conjugate as a complex octonion in 
$\OO\otimes_\R\C$) i.e. for every ${\bf v}=v_1\ee_1+\dots+v_7\ee_7$ we define 
$\overline{\bf v}:=\overline{v}_1\ee_1+\dots+\overline{v}_7\ee_7$. 

Moreover $O'(\Lambda)$ is homeomorphic to $S^6$, does not depend on 
the Samelson parameter $u=a+\sqrt{-1}b$ with $a\not=0$ and in particular 
$f^\C(-\ee_1)f^\C(-\ee_1)=\Lambda\in O'(\Lambda)$ justifying the notation. 
\label{komplexegyenloseg}
\end{lemma}

\begin{proof}
We have seen that $\Lambda^2=\Lambda^{-1}=f(\ee_1)^{-1}=f(-\ee_1)=
f^\C(-\ee_1)$ hence knowing that $s_1\in S_u\subset{\rm G}_2^\C$ is an 
automorphism we find 
$f^\C({\bf z})=s_1\Lambda^2s_1^{-1}=s_1f^\C(-\ee_1)s_1^{-1}=f^\C(s_1(-\ee_1))$ 
hence ${\bf z}=-s_1\ee_1$. Consider the unique decomposition $s_1=\lim\limits_n
\big(A^{1/n}h_1^{1/n}\big)^n$ where $h_1\in S_u\cap{\rm SU}(3)^\C$ and 
using the previously constructed root basis 
for ${\mathfrak s}_u$ we introduced $A:=
\exp\big(\sqrt{-6}\:z_2U_{+1}-\sqrt{-6}\:z_4U_{-2}+\sqrt{6}\:z_6U_{+3}\big)
\in S_u$. This latter matrix is relatively easy to compute because fortunately 
$U_{\pm k}^3=0$ for $k=1,2,3$ yielding 
\[A=\left(\begin{smallmatrix}
   1 & z_2 &-\sqrt{-1}z_2 & z_4 & -\sqrt{-1}z_4 & z_6 & -\sqrt{-1}z_6\\
   -z_2 & 1-\frac{1}{2}z_2^2 & \frac{\sqrt{-1}}{2}z_2^2 & 
\frac{\sqrt{-1}}{2}z_6-\frac{1}{2}z_2z_4 & \frac{1}{2}z_6+\frac{\sqrt{-1}}{2}
z_2z_4 & \frac{\sqrt{-1}}{2}z_4-\frac{1}{4}z_2z_6 & 
\frac{1}{2}z_4+\frac{\sqrt{-1}}{4}z_2z_6\\
\sqrt{-1}z_2 & \frac{\sqrt{-1}}{2}z_2^2 & 1+\frac{1}{2}z_2^2 &  
\frac{1}{2}z_6+\frac{\sqrt{-1}}{2}z_2z_4 & -\frac{\sqrt{-1}}{2}z_6+
\frac{1}{2}z_2z_4 & -\frac{1}{2}z_4+\frac{\sqrt{-1}}{4}z_2z_6 &  
\frac{\sqrt{-1}}{2}z_4+\frac{1}{4}z_2z_6 \\
-z_4 & -\frac{\sqrt{-1}}{2}z_6-\frac{1}{2}z_2z_4 & 
-\frac{1}{2}z_6+\frac{\sqrt{-1}}{2}z_2z_4 & 1-\frac{1}{2}z_4^2 & 
\frac{\sqrt{-1}}{2}z_4^2 & -\frac{\sqrt{-1}}{2}z_2-\frac{1}{4}z_4z_6
&  -\frac{1}{2}z_2+\frac{\sqrt{-1}}{4}z_4z_6\\
\sqrt{-1}z_4 & -\frac{1}{2}z_6+\frac{\sqrt{-1}}{2}z_2z_4 & 
\frac{\sqrt{-1}}{2}z_6+\frac{1}{2}z_2z_4 & 
\frac{\sqrt{-1}}{2}z_4^2 & 1+\frac{1}{2}z_4^2 & 
\frac{1}{2}z_2+\frac{\sqrt{-1}}{4}z_4z_6 & 
-\frac{\sqrt{-1}}{2}z_2+\frac{1}{4}z_4z_6 \\ 
-z_6 & -\frac{\sqrt{-1}}{2}z_4-\frac{1}{4}z_2z_6 & 
\frac{1}{2}z_4+\frac{\sqrt{-1}}{4}z_2z_6 & \frac{\sqrt{-1}}{2}z_2-
\frac{1}{4}z_4z_6 & -\frac{1}{2}z_2+\frac{\sqrt{-1}}{4}z_4z_6 &  
1-\frac{1}{2}z_6^2 & \frac{\sqrt{-1}}{2}z_6^2\\
\sqrt{-1}z_6 & -\frac{1}{2}z_4+\frac{\sqrt{-1}}{4}z_2z_6 & 
-\frac{\sqrt{-1}}{2}z_4+\frac{1}{4}z_2z_6 & 
\frac{1}{2}z_2+\frac{\sqrt{-1}}{4}z_4z_6 & 
\frac{\sqrt{-1}}{2}z_2+\frac{1}{4}z_4z_6 & 
\frac{\sqrt{-1}}{2}z_6^2 & 1+\frac{1}{2}z_6^2
\end{smallmatrix}\right)\]
thus it is algebraic. Since $h_1$ is the stabilizer of $\ee_1$ and 
$\dim_\C S_u\Lambda^2S_u^{-1}=3$ as well as $A$ already depends on $3$ 
complex parameters by setting $h_1:=e$ we can replace $s_1$ with the simple 
matrix $A$ at the price of parameterizing an open subset of 
$S_u\Lambda^2S_u^{-1}$ only. Thus the action of this $A$ on $-\ee_1$ gives
\[{\bf z}=-{\bf e}_1+z_2{\bf e}_2-\sqrt{-1}z_2{\bf e}_3+z_4{\bf e}_4
          -\sqrt{-1}z_4{\bf e}_5+z_6{\bf e}_6-\sqrt{-1}z_6{\bf e}_7\]
thus ${\bf z}\in (S^6)^\C$ with arbitrary parameters $z_2,z_4,z_6\in\C$. 
Likewise $f^\C({\bf v})\in\overline{S}_u\Lambda^2\overline{S}_u^{-1}$ if  
\[{\bf v}=-{\bf e}_1+v_2{\bf e}_2+\sqrt{-1}v_2{\bf e}_3+v_4{\bf e}_4
          +\sqrt{-1}v_4{\bf e}_5+v_6{\bf e}_6+\sqrt{-1}v_6{\bf e}_7\]
hence ${\bf v}\in (S^6)^\C$ with further $v_2,v_4,v_6\in\C$. In this case for 
the element $s_2\in\overline{S}_u$ taking $-\ee_1$ to ${\bf v}$ it is enough 
to put $B:=\exp\big(-\sqrt{-6}\:v_2U_{-1}+\sqrt{-6}\:v_4U_{+2}+
\sqrt{6}\:v_6U_{-3}\big)\in\overline{S}_u$ which looks like 
\[B=\left(\begin{smallmatrix}
1 & v_2 & \sqrt{-1}v_2 & v_4 & \sqrt{-1}v_4 & v_6 & \sqrt{-1}v_6 \\
-v_2 & 1-\frac{1}{2}v_2^2 & -\frac{\sqrt{-1}}{2}v_2^2 & 
-\frac{\sqrt{-1}}{2}v_6-\frac{1}{2}v_2v_4 & \frac{1}{2}v_6-
\frac{\sqrt{-1}}{2}v_2v_4 & -\frac{\sqrt{-1}}{2}v_4-\frac{1}{4}v_2v_6 &
\frac{1}{2}v_4-\frac{\sqrt{-1}}{4}v_2v_6 \\ 
-\sqrt{-1}v_2 & -\frac{\sqrt{-1}}{2}v_2^2 & 1+\frac{1}{2}v_2^2 & 
\frac{1}{2}v_6-\frac{\sqrt{-1}}{2}v_2v_4 & \frac{\sqrt{-1}}{2}v_6+
\frac{1}{2}v_2v_4 & -\frac{1}{2}v_4-\frac{\sqrt{-1}}{4}v_2v_6 & 
-\frac{\sqrt{-1}}{2}v_4+\frac{1}{4}v_2v_6 \\ 
-v_4 & \frac{\sqrt{-1}}{2}v_6-\frac{1}{2}v_2v_4 & 
-\frac{1}{2}v_6-\frac{\sqrt{-1}}{2}v_2v_4 & 1-\frac{1}{2}v_4^2 & 
-\frac{\sqrt{-1}}{2}v_4^2 & \frac{\sqrt{-1}}{2}v_2-\frac{1}{4}v_4v_6 & 
-\frac{1}{2}v_2-\frac{\sqrt{-1}}{4}v_4v_6 \\
-\sqrt{-1}v_4 & -\frac{1}{2}v_6-\frac{\sqrt{-1}}{2}v_2v_4 
& -\frac{\sqrt{-1}}{2}v_6+\frac{1}{2}v_2v_4 & -\frac{\sqrt{-1}}{2}v_4^2 & 
1+\frac{1}{2}v_4^2 & \frac{1}{2}v_2-\frac{\sqrt{-1}}{4}v_4v_6 & 
\frac{\sqrt{-1}}{2}v_2+\frac{1}{4}v_4v_6 \\ 
-v_6 & \frac{\sqrt{-1}}{2}v_4-\frac{1}{4}v_2v_6 & \frac{1}{2}v_4-
\frac{\sqrt{-1}}{4}v_2v_6 & -\frac{\sqrt{-1}}{2}v_2-\frac{1}{4}v_4v_6 & 
-\frac{1}{2}v_2-\frac{\sqrt{-1}}{4}v_4v_6 & 1-\frac{1}{2}v_6^2 & 
-\frac{\sqrt{-1}}{2}v_6^2 \\ 
-\sqrt{-1}v_6 & -\frac{1}{2}v_4-\frac{\sqrt{-1}}{4}v_2v_6 & 
\frac{\sqrt{-1}}{2}v_4+\frac{1}{4}v_2v_6 & \frac{1}{2}v_2-
\frac{\sqrt{-1}}{4}v_4v_6 & -\frac{\sqrt{-1}}{2}v_2+\frac{1}{4}v_4v_6 & 
-\frac{\sqrt{-1}}{2}v_6^2 & 1+\frac{1}{2}v_6^2
\end{smallmatrix}\right)\]
akin to $A$. The perturbed complexified orbit 
$O'(\Lambda)^\C$ as defined in (\ref{komplexpalya'}) is of the form 
$S_u\Lambda^2S_u^{-1}\overline{S}_u\Lambda^2\overline{S}_u^{-1}$ 
consequently we obtain an open subset of it which is however obviously 
closed too; therefore we find that 
$O'(\Lambda)^\C=\{f^\C({\bf z})f^\C({\bf v})\:\vert\:
\mbox{${\bf z},{\bf v}\in (S^6)^\C$ as above}\}$. 

Let us turn now to $O'(\Lambda)$. As defined in (\ref{palya'}) the 
``real part'' $O'(\Lambda)$ arises within $O'(\Lambda)^\C$ by imposing on the 
pairs $s_1\in S_u$ and $s_2\in\overline{S}_u$ the further reality 
condition $s_1=s=\overline{s}_2$. 
This is because we know that 
$g=\lim\limits_n\big(s^{1/n}\overline{s}^{1/n}\big)^n$ is a decomposition of 
a real element; we also know that given $g\in {\rm G}_2$ one can find 
an element $s\in S_u$ such that for every $h\in {\rm SU}(3)^\C$ 
there exist $h_1\in S_u\cap{\rm SU}(3)^\C$ and 
$h_2\in\overline{S}_u\cap{\rm SU}(3)^\C$ satisfying  
$gh=\lim\limits_n\big(s^{1/n}h_1^{1/n}\overline{s}^{1/n}h_2^{1/n}\big)^n$. 
Thus as $s$ runs over $S_u$ the correspondence 
\[s\Lambda^2s^{-1}\:\overline{s}\Lambda^2\overline{s}^{-1}
\:\:\Longleftrightarrow\:\:
\lim\limits_n\big(s^{1/n}\:\overline{s}^{1/n}\big)^n\:\Lambda\:
\lim\limits_n\big(s^{1/n}\:\overline{s}^{1/n}\big)^{-n}\]
is a homeomorphism between $O'(\Lambda)$ as defined in (\ref{palya'}) 
and $O(\Lambda)$ as defined in (\ref{palya}) i.e. between $O'(\Lambda)$ and 
$S^6$. To find an explicit parameterization of $O'(\Lambda)$ we go on again 
with observing that every element $g\in{\rm G}_2$ can be written 
non-uniquely in the form 
$g=\lim\limits_n\big(s^{1/n}\overline{s}^{1/n}\big)^n$ with $s\in S_u$; 
and in addition, every $s\in S_u$ admits a unique decomposition 
$s=\lim\limits_n\big(A^{1/n}h^{1/n}\big)^n$ where $A$ is the matrix above and 
$h\in S_u\cap{\rm SU}(3)^\C$. Introducing 
$Z(z_2,z_4,z_6):=\sqrt{-6}\:z_2U_{+1}-\sqrt{-6}\:z_4U_{-2}+\sqrt{6}\:z_6U_{+3}
\in{\mathfrak g}_2^\C$ hence $\exp Z(z_2,z_4,z_6)=A(z_2,z_4,z_6)$ 
as so far and likewise $W\in{\mathfrak g}_2^\C$ via $\exp W=h$ these Trotter 
decompositions imply that $g\in {\rm G}_2$ can be written as
\[g=\exp\big(Z(z_2,z_4,z_6)+W+\overline{Z}(\overline{z}_2,
\overline{z}_4,\overline{z}_6)+\overline{W}\big)\:\:.\]
Consequently, taking into account the connectedness and compactness of 
${\rm G}_2$ and that the injective restriction of 
$\exp:{\mathfrak g}_2\rightarrow{\rm G}_2$ is a proper map, 
since $A(0,0,0)=e$ there exists a smallest bounded 
connected open subset $\Omega\subset\C^3$ about the origin, 
the ``parameter space'', such that every element of ${\rm G}_2$ 
can be obtained out of a matrix $A(z_2,z_4,z_6)$ satisfying $(z_2,z_4,z_6)\in
\overline{\Omega}$ and an element $h$ also belonging to the closure of a 
neighbourhood of $e\in S_u\cap{\rm SU}(3)^\C$. Since 
$\dim_\R O'(\Lambda)=6$ and the matrix $A(z_2,z_4,z_6)$ already depends on $6$ 
real variables, we can represent an open subset of $O'(\Lambda)$ 
simply by matrices of the form 
$A\Lambda^2A^{-1}\overline{A}\Lambda^2\overline{A}^{-1}$. 
It is easy to see that this subset is in fact closed too thus coincides with 
$O'(\Lambda)$. Because $f^\C({\bf z})
=A\Lambda^2A^{-1}$ hence $f^\C(\overline{\bf z})=\overline{f^\C({\bf z})}
=\overline{A\Lambda^2A^{-1}}=\overline{A}\Lambda^2\overline{A}^{-1}$ 
we find that the claimed description $O'(\Lambda)
=\{f^\C({\bf z})f^\C(\overline{\bf z})\:\vert\:\mbox{with 
${\bf z}\in (S^6)^\C$ as above with $(z_2,z_4,z_6)\in\overline{\Omega}$}\}$ 
follows. 

Finally observe that $O'(\Lambda)$ is homeomorphic to $S^6$ as we already know 
and it does not depend on $u=a+\sqrt{-1}b$. Moreover 
$z_{2}=z_4=z_6=0$ gives ${\bf z}=-\ee_1$ and $f^\C(-\ee_1)f^\C(-\ee_1)=
\Lambda^2\Lambda^2=\Lambda\in O'(\Lambda)$ as stated. 
\end{proof}
  

\noindent Our last steps are then as follows. Put 
\[F({\bf z},\overline{\bf z}):=f^\C({\bf z})f^\C(\overline{\bf z})\]
and let $\{X_1,\dots,X_6\}$ be a reasonable $\R$-basis at 
$T_{F({\bf z},\overline{\bf z})}O'(\Lambda)$ and extend it smoothly to an 
$\R$-frame over the punctured space $O'(\Lambda)\setminus\{{\rm point}\}$ 
permitted by the parameterization of Lemma \ref{komplexegyenloseg} i.e. take 
all ${\bf z}=-{\bf e}_1+z_2{\bf e}_2-\sqrt{-1}z_2{\bf e}_3+z_4{\bf e}_4
          -\sqrt{-1}z_4{\bf e}_5+z_6{\bf e}_6-\sqrt{-1}z_6{\bf e}_7$ 
with $(z_2,z_4,z_6)\in\Omega$ (note that no frame can extend further 
because $TO'(\Lambda)\cong TS^6$ is not trivial). Then the construction of 
$J_u$ at ${\mathfrak g}_2$ allows one to compute the action 
of $J_u$ at $T_{F({\bf z},\overline{\bf z})}O'(\Lambda)$ by left-invariance 
i.e. applying the formula $J_uX_i:=L_{F({\bf z},\overline{\bf z})*}J_u
L_{F({\bf z},\overline{\bf z})*}^{-1}X_i$ or in simpler notation 
\[J_uX_i:=F({\bf z},\overline{\bf z})J_u F({\bf z},\overline{\bf z})^{-1}X_i
\:\:.\] 
To carry out these computations one has to proceed in principle as follows: 
first by the aid of the Hermitian scalar product 
$\langle\:\cdot\:,\:\cdot\:\rangle^\C$ on ${\mathfrak g}_2^\C$ one computes 
the matrix coefficients 
$\big\langle F({\bf z},\overline{\bf z})^{-1}X_i\:,\:H_{\pm}
\big\rangle^\C\in\C$, etc. in order to expand the pullbacks 
$F({\bf z},\overline{\bf z})^{-1}X_i\in{\mathfrak g}_2^\C$ in the $\C$-linear
extension of the ${\mathfrak g}_2$-basis 
$\{H_{\pm}, X_{\pm k},Y_{\pm k}\}$ used in (\ref{Jhatas}); 
secondly one uses the explicit action (\ref{Jhatas}) of $J_u$ in this basis 
to obtain $J_uF({\bf z},\overline{\bf z})^{-1}X_i\in{\mathfrak g}_2^\C$; and 
thirdly by multiplying this matrix with $F({\bf z},\overline{\bf z})$ from 
the left one transfers the result from ${\mathfrak g}_2^\C$ back to 
$T_{F({\bf z},\overline{\bf z})}{\rm G}_2^\C$. Writing 
\[J_{ij}:=\big\langle J_uF({\bf z},\overline{\bf z})^{-1}X_i\:,\:
F({\bf z},\overline{\bf z})^{-1}X_j\big\rangle^\C\]
one tries to check that $J_uX_i-\sum\limits_{j=1}^6J_{ij}X_j=0$ implying 
$F({\bf z},\overline{\bf z})J_u
F({\bf z},\overline{\bf z})^{-1}X_i\in T_{F({\bf z},\overline{\bf z})}
O'(\Lambda)\subset T_{F({\bf z},\overline{\bf z})}{\rm G}_2^\C$ i.e. the 
Samelson almost complex structure indeed restricts to 
$O'(\Lambda)\setminus\{{\rm point}\}$. Moreover introducing the 
singular metric 
\[g_{ik}:=\big\langle F({\bf z},\overline{\bf z})^{-1}X_i\:,\:
F({\bf z},\overline{\bf z})^{-1}X_k\big\rangle^\C\]
and its inverse $g^{ik}$ along $O'(\Lambda)\setminus\{{\rm point}\}$ 
one tries to check that the matrix coefficients 
\[J^i_k:=\sum\limits_{j=1}^6g^{ij}J_{jk}\]
of the true $(1,1)$-tensor $J_u$ are already smooth as 
$(z_2,z_4,z_6)\rightarrow\partial\overline{\Omega}$  
yielding a well-defined $J_u$ over the whole $O'(\Lambda)$ as desired.
However it is a {\it mission impossible} to perform and plot these brute force 
computations here.
 
The Samelson complex structure and its integrability shows up more clearly if 
we rather pass to complexification and exploit 
some technical observations made during the proof of Theorem \ref{letezes}. 
The obvious tangent vectors $\Big\{\frac{\partial 
F({\bf z},\overline{\bf z})}{\partial z_{2}},\frac{\partial 
F({\bf z},\overline{\bf z})}{\partial \overline{z}_{2}},\dots,
\frac{\partial F({\bf z},\overline{\bf z})}{\partial z_{6}},\frac{\partial
F({\bf z},\overline{\bf z})}{\partial \overline{z}_{6}}\Big\}$ 
comprise a $\C$-basis in $T_{F({\bf z},\overline{\bf z})}O'(\Lambda)^\C$. 
Consequently, writing $L_{F({\bf z},\overline{\bf z})*}^{-1}\frac{\partial
F({\bf z},\overline{\bf z})}{\partial z_{2k}}$ etc., simply as a 
matrix product $F({\bf z},\overline{\bf z})^{-1}\frac{\partial
F({\bf z},\overline{\bf z})}{\partial z_{2k}}$ etc., from now on, 
$\Big\{F({\bf z},\overline{\bf z})^{-1}\frac{\partial
F({\bf z},\overline{\bf z})}{\partial z_{2}},\dots,F({\bf z},
\overline{\bf z})^{-1}\frac{\partial
F({\bf z},\overline{\bf z})}{\partial \overline{z}_{6}}\Big\}$ is a 
$\C$-basis in $L_{F({\bf z},\overline{\bf z})*}^{-1}T_{F({\bf z},
\overline{\bf z})}O'(\Lambda)^\C\subset{\mathfrak g}_2^\C$. It has an 
asymmetric shape 
\[\left\{f^\C(\overline{\bf z})^{-1}f^\C({\bf z})^{-1}
\frac{\partial f^\C({\bf z})}{\partial z_{2}}
f^\C(\overline{\bf z})\:,\:
f^\C(\overline{\bf z})^{-1}\frac{\partial
f^\C(\overline{\bf z})}{\partial\overline{z}_{2}}\:,\dots,
f^\C(\overline{\bf z})^{-1}\frac{\partial
f^\C(\overline{\bf z})}{\partial\overline{z}_{6}}\right\}\:\:.\] 
Instead of this straightforward basis consider the not only symmetric 
but even simpler collection 
\[\left\{f^\C({\bf z})^{-1}
\frac{\partial f^\C({\bf z})}{\partial z_{2}}\:,\:
f^\C(\overline{\bf z})^{-1}\frac{\partial
f^\C(\overline{\bf z})}{\partial\overline{z}_{2}}\:,\dots,
f^\C({\bf z})^{-1}\frac{\partial f^\C({\bf z})}{\partial z_{6}}\:,\:
f^\C(\overline{\bf z})^{-1}\frac{\partial
f^\C(\overline{\bf z})}{\partial\overline{z}_{6}}\right\}\]
within ${\mathfrak g}_2^\C$ (but yet the full size of these matrices can be 
guessed from the shape of the real derivative $f({\bf x})^{-1}f({\bf x})_*\xi$ 
computed above). Note that actually $\Big\{f^\C({\bf z})^{-1}
\frac{\partial f^\C({\bf z})}{\partial z_{2}}\:,\:f^\C({\bf z})^{-1}
\frac{\partial f^\C({\bf z})}{\partial z_{4}}\:,\:f^\C({\bf z})^{-1}
\frac{\partial f^\C({\bf z})}{\partial z_{6}}\Big\}$ gives a $\C$-basis in 
$L_{f^\C({\bf z})*}^{-1}T_{f^\C({\bf z})}
S_u\Lambda^2S_u^{-1}\subset{\mathfrak s}_u$ and the rest forms a $\C$-basis in 
$L_{f^\C(\overline{\bf z})*}^{-1}T_{f^\C(\overline{\bf z})}
\overline{S}_u\Lambda^2\overline{S}_u^{-1}\subset\overline{\mathfrak s}_u$. 
Consequently, since by recalling (\ref{hasitas}) we know that a splitting 
\[L_{F({\bf z},\overline{\bf z})*}^{-1}T_{F({\bf z},
\overline{\bf z})}O'(\Lambda)^\C=L_{f^\C({\bf z})*}^{-1}T_{f^\C({\bf z})}
S_u\Lambda^2S_u^{-1}\oplus L_{f^\C(\overline{\bf z})*}^{-1}
T_{f^\C(\overline{\bf z})}\overline{S}_u\Lambda^2\overline{S}_u^{-1}\]
holds, in fact this symmetric collection of matrices constitutes a $\C$-basis 
in $L_{F({\bf z},\overline{\bf z})*}^{-1}T_{F({\bf z},
\overline{\bf z})}O'(\Lambda)^\C$ too. Define $Z_{2k-1}:=
F({\bf z},\overline{\bf z})f^\C({\bf z})^{-1}
\frac{\partial f^\C({\bf z})}{\partial z_{2k}}$ and $Z_{2k}:=
F({\bf z},\overline{\bf z})f^\C(\overline{\bf z})^{-1}\frac{\partial
f^\C(\overline{\bf z})}{\partial\overline{z}_{2k}}$ for $k=1,2,3$; 
then $\{Z_1,\dots, Z_6\}$ gives rise to a smooth $\C$-frame field if
$(z_2,z_4,z_6)\in\Omega$ that is over 
the punctured space $O'(\Lambda)\setminus\{{\rm point}\}$.

From this point we proceed as before. Introducing again 
$J^\C_{ij}:=\big\langle J^\C_uF({\bf z},\overline{\bf z})^{-1}Z_i\:,\:
F({\bf z},\overline{\bf z})^{-1}Z_j\big\rangle^\C$ we simply find that 
\[J^\C_{ij}=\pm\sqrt{-1}\:g_{ij}\] 
where we write again $g_{kl}:=\big\langle F({\bf z},\overline{\bf z})^{-1}Z_k
\:,\:F({\bf z},\overline{\bf z})^{-1}Z_l\big\rangle^\C$
and $g^{kl}$ for the components of a complex singular metric and its inverse
along $O'(\Lambda)\setminus\{{\rm point}\}$. This is because from the 
observations on the frame made above we know that 
$F({\bf z},\overline{\bf z})^{-1}Z_k\in{\mathfrak s}_u={\mathfrak g}_2^{1,0}$ 
for $k=1,3,5$ and likewise $F({\bf z},\overline{\bf z})^{-1}Z_k
\in\overline{\mathfrak s}_u={\mathfrak g}_2^{0,1}$ for $k=2,4,6$ where 
${\mathfrak g}_2^\C={\mathfrak g}_2^{1,0}\oplus{\mathfrak g}_2^{0,1}$ is the 
$\pm\sqrt{-1}$-eigenspace decomposition of the complexified Lie algebra with 
respect to the complex linear extension $J^\C_u$ 
of the Samelson complex structure. Consequently by left-invariance of the 
Samelson complex structure the matrix elements  
${J^\C}^i_k=\sum\limits_{j=1}^6g^{ij}J^\C_{jk}$ of the true 
$(1,1)$-tensor $J^\C_u$ simply look like 
\[{J^\C}^i_k=\pm\sqrt{-1}\:\delta^i_k\] 
hence it readily follows that despite the degeneration of the 
frame $\{Z_1,\dots, Z_6\}$ the $(1,1)$-tensor $J^\C_u$ itself remains smooth as 
$(z_2,z_4,z_6)\rightarrow\partial\overline{\Omega}$ hence 
we obtain a well-defined $J^\C_u$ over the whole $O'(\Lambda)$. 

Reality of $O'(\Lambda)\subset O'(\Lambda)^\C$ implies that 
$F({\bf z},\overline{\bf z})
\overline{F({\bf z},\overline{\bf z})^{-1}Z}_{2k-1}=F({\bf z},\overline{\bf z})
F({\bf z},\overline{\bf z})^{-1}Z_{2k}$ for $k=1,2,3$ 
in $T_{F({\bf z},\overline{\bf z})}O'(\Lambda)^\C=
T_{F({\bf z},\overline{\bf z})}O'(\Lambda)\otimes_\R\C$. 
Consequently we can introduce an $\R$-frame over the punctured space 
$O'(\Lambda)\setminus\{{\rm point}\}$ canonically induced by the $\C$-frame 
$\{Z_1,\dots,Z_6\}$ i.e. $X_{2k-1}:=\frac{1}{2}(Z_{2k-1}+Z_{2k})$ and similarly 
$X_{2k}:=\frac{1}{2\sqrt{-1}}(Z_{2k-1}-Z_{2k})$ for all $k=1,2,3$. In the 
resulting $\R$-frame $\{X_1,\dots,X_6\}$ the real almost 
complex tensor field therefore takes the very simple form  
\[J_u\vert_{O'(\Lambda)\setminus\{{\rm point}\}}=\left(\begin{matrix}
                     0 & 1 & 0 & 0 & 0 & 0\\
                     -1 & 0 & 0 & 0 & 0 & 0\\
                     0 & 0 & 0 & 1 & 0 & 0\\ 
                     0 & 0 & -1 & 0 & 0 & 0\\
                     0 & 0 & 0 & 0 & 0 & 1\\
                     0 & 0 & 0 & 0 & -1 & 0\\
                  \end{matrix}\right)\] 
demonstrating its smooth extendibility over the 
whole $O'(\Lambda)\cong S^6$. Concerning its integrability however 
some care is needed because the frame $\{X_1,\dots,X_6\}$ is not obviously 
torsion-free hence the constancy of $J_u$ in this frame might not prove 
anything. Let us therefore compute the Nijenhuis tensor. Observe 
that the $(0,1)$-parts $X_j^{1,0}=\frac{1}{2}(X_j-\sqrt{-1}J_uX_j)$ look 
simply $X^{1,0}_{2k-1}=\frac{1}{2}Z_{2k-1}$ and 
$X^{1,0}_{2k}=\frac{1}{2\sqrt{-1}}Z_{2k-1}$ for $k=1,2,3$ consequently 
$F({\bf z},\overline{\bf z})^{-1}X_j^{1,0}\in{\mathfrak g}_2^{1,0}$ for 
all $j=1,\dots,6$. Thus 
\[N_{J_u}(X_i,X_j)=\big[X_i^{1,0}\:,\:X_j^{1,0}\big]^{0,1}=0\]
because $\big[{\mathfrak g}_2^{1,0}\:,{\mathfrak g}_2^{1,0}\big]^{0,1}=0$. 
Thus $N_{J_u}=0$ guaranteeing integrability over 
$O'(\Lambda)\setminus\{{\rm point}\}$ hence over the whole $O'(\Lambda)$.
 
Note that in this second picture all the 
computational difficulties and geometric subtleties concerning $J_u$ 
have been compressed into the highly non-trivial frame field 
$\{X_1,\dots,X_6\}$ in which $J_u$ looks simply constant. For instance, since 
as exhibited in Lemma \ref{komplexegyenloseg} the space $O'(\Lambda)$ hence 
the frame $\{X_1,\dots,X_6\}$ along $O'(\Lambda)\setminus\{{\rm point}\}$ are 
independent of the Samelson moduli parameter $u=a+\sqrt{-1}b$ with $a\not=0$, 
we see that the constructed complex manifold $X_u$ in Theorem \ref{letezes} is 
{\it unique} i.e. independent of the Samelson moduli parameter 
$u\in\C\setminus\sqrt{-1}\R$. 

With this observation we conclude the struggle with the explicit 
construction of the complex structure on the six-sphere.


\section{Appendix: inner automorphisms and $\pi_6({\rm G}_2)\cong\Z_3$}
\label{four}


To close we take a look of the conjugate orbit 
$O(\Lambda )\subset{\rm G}_2$ 
in (\ref{palya}) from a topological viewpoint. As we have seen in 
Lemma \ref{egyenloseg} it can be identified with the image of the map 
$f:S^6\rightarrow {\rm G}_2$ in (\ref{belso.automorfizmusok}); we 
demonstrate that in this form the conjugate orbit represents the generator of 
the sixth homotopy group of the automorphism group of the octonions. 
Consequently this homotopy group is non-trivial and is generated by inner 
automorphisms (thus $O'(\Lambda)\subset{\rm G}_2^\C$ in (\ref{palya'}) 
is also homotopically non-trivial). We acknowledge that this group has been 
known for a long time \cite{mim} and even our proof is 
essentially the same as the nice geometric one in \cite{cha-rig}. 

\begin{theorem}{\rm (cf. \cite{cha-rig, mim})} 
There exists an isomorphism $\pi_6({\rm G}_2)\cong\Z_3$\:. Moreover the 
map (\ref{belso.automorfizmusok}) constructed out of the collection of 
rotations induced by inner automorphisms (\ref{konjugalas}) of the 
octonions, is a representative of the generator of this homotopy group. 
\label{homotopia}
\end{theorem} 

\begin{proof}
Recall that ${\rm Spin}(7)\subset{\rm Cliff}_0(\R^7)\cong{\rm Cliff}(\R^6)
\cong\R (8)$ hence the unique spin representation of ${\rm SO}(7)$ acts on 
$\R^8$. This gives rise to an embedding ${\rm Spin}(7)\subset{\rm SO}(8)$. 
The projection $P:{\rm SO}(8)\rightarrow S^7$ sending a matrix onto 
its (let us say) first column restricts to ${\rm Spin}(7)$ providing us 
with a projection $\tilde{p}:{\rm Spin}(7)\rightarrow S^7$. Dividing 
this with the center $\Z_2\cong Z({\rm SO}(8))\subset{\rm SO}(8)$ 
we obtain another projection $p:{\rm SO}(7)\rightarrow\R P^7$. The 
geometric meaning of this map is straightforward: the preimage of a point of 
$\R P^7$ i.e. a line in $\R^8$ consists of those rotations which keep 
this line fixed therefore act only on a hyperplane perpendicular to this 
line: dimension counting shows that these transformations are exactly the 
automorphisms of the octonions hence their collection is isomorphic to 
${\rm G}_2$. Consequently the projection $p:{\rm SO}(7)\rightarrow\R P^7$ 
is the classical ${\rm G}_2$-fibration of ${\rm SO}(7)$ over $\R P^7$. It has 
an associated homotopy exact sequence whose relevant segment for us is  
\[\xymatrix{& 
& \pi_7({\rm SO}(7),{\rm G}_2,e)\ar[d]_{p_*}^{\cong}\ar[dr]^{\partial_*}& & &\\ 
\dots\ar[r]&\pi_7({\rm SO}(7),e)\ar[ur]^{j_*}&\pi_7(\R P^7,p(e))&
\pi_6({\rm G}_2,e)\ar[r]&\pi_6({\rm SO}(7),e)\ar[r]&\dots}\]
where $j:({\rm SO}(7),e,e)\rightarrow ({\rm SO}(7),{\rm G}_2,e)$ is 
induced  by the embedding $e\in{\rm G}_2\subset{\rm SO}(7)$. 

First let us compute $\pi_6({\rm G}_2)$. Take $S^7=\{2\cos t\cdot\ee_0+2\sin t
\cdot {\bf x}\:\vert\:0\leqq t\leqq\pi\:,\: {\bf x}\in S^6\}$. The conjugation 
(\ref{konjugalas}) can be enhanced to an orthogonal transformation of the 
octonions which on a particular ${\bf u}\in\OO$ has the form (again parantheses 
omitted) 
\begin{equation} 
{\bf u}\mapsto (2\cos t\cdot\ee_0 +2\sin t\cdot{\bf x}) 
{\bf u}(2\cos t\cdot\ee_0 + 2\sin t\cdot {\bf x})^{-1}=
(\cos t\cdot\ee_0+\sin t\cdot {\bf x})
{\bf u}(\cos t\cdot\ee_0-\sin t\cdot {\bf x})
\label{nagykonjugalas} 
\end{equation} 
and provides us with a map from $S^7$ into ${\rm SO}(8)$. However this 
apparent ${\rm SO}(8)$ transformation of $\OO\cong\R^8$ leaves $\re
\OO\cong\R$ invariant i.e. acts only on $\im\OO\cong\R^7$ 
therefore it is actually an ${\rm SO}(7)$ transformation. This way we 
obtain a map $F :S^7\rightarrow{\rm SO}(7)$ such that $[F]=1\in\pi_7({\rm 
SO}(7))\cong\Z$ i.e. its homotopy class is a generator 
\cite{tod-sai-yok}. Take now the $7$-cell $e^7:=\{2\cos t\cdot\ee_0+2\sin t
\cdot {\bf x}\:\vert\:0\leqq t\leqq\frac{\pi}{3}\:,\:{\bf x}\in S^6\}
\subset S^7$. Its boundary is $\partial e^7=\{\ee_0+\sqrt{3}\:
{\bf x}\:\vert\:{\bf x}\in S^6\}$ 
hence constitutes the inner automorphisms (\ref{konjugalas}) of the 
octonions therefore via (\ref{nagykonjugalas}) it lies within ${\rm G}_2 
\subset{\rm SO}(7)$. Consequently restriction to this $7$-cell gives rise to a 
map $F\vert_{e^7} :(e^7,\partial e^7)\rightarrow ({\rm SO}(7),{\rm G}_2)$ 
satisfying $[F\vert_{e^7}]=1\in\pi_7({\rm SO}(7),{\rm G}_2)\cong\pi_7(\R 
P^7)\cong\Z$ i.e. its homotopy class continues to be a generator. Taking 
into account that the third power of an inner automorphism is the identity 
it is clear that $j_*[F]=3[F\vert_{e^7}]$. Consequently, since $\pi_6({\rm 
SO}(7))\cong 0$ we conclude from the homotopy exact sequence that 
$\pi_6({\rm G}_2)\cong\Z/3\Z$, written as $\Z_3$, as desired.

Regarding the generator, it readily follows from (\ref{konjugalas}) and 
(\ref{nagykonjugalas}) that $\partial (F\vert_{e^7})=f$ where 
$f:S^6\rightarrow {\rm G}_2$ is the map (\ref{belso.automorfizmusok}). 
Therefore $\partial_*[F\vert_{e^7}]=[f]\in\pi_6({\rm G}_2)$. By exactness 
$\partial_*\not=0$ and $F\vert_{e^7}$ represents the generator, hence its 
image $f$ also represents a non-trivial element in $\pi_6({\rm G}_2)$ which is 
the generator.
\end{proof}

\end{document}